\documentclass[a4paper,10pt]{article}
\title{The Dixmier map for nilpotent super Lie algebras
}
\author{Estanislao Herscovich
\thanks{The author is an Alexander von Humboldt fellow.}}
\date{}

\usepackage{amsmath,amsthm,amsfonts,amssymb,stmaryrd}
\usepackage{latexsym}
\usepackage{color}
\usepackage{graphicx}
\usepackage{float}  
\usepackage{fullpage}
\usepackage[small, bf, margin=90pt, tableposition=bottom]{caption}
\usepackage[T1]{fontenc}
\usepackage{palatino}
\usepackage[fleqn,tbtags]{mathtools}

\usepackage[
  breaklinks,
  naturalnames,
  ]{hyperref}

\usepackage[alphabetic]{amsrefs}

\input{xy}
\xyoption{all}
\xyoption{poly}
\usepackage[all]{xy}

\newtheorem{teo}{Theorem}[section]
\newtheorem{teo*}{Theorem}
\newtheorem{coro}[teo]{Corollary}
\newtheorem{lema}[teo]{Lemma}
\newtheorem{prop}[teo]{Proposition}
\newtheorem{rem}[teo]{Remark}

\newtheorem{ex}[teo]{Example}

\numberwithin{equation}{section}                    

\def\id{{\mathrm{id}}}
\let\oldqed\qed
\renewcommand\qed{\oldqed\par\bigskip}

\newcommand\cl[1]{{\langle#1\rangle}}

\newcommand\ZZ{{\mathbb{Z}}}
\newcommand\NN{{\mathbb{N}}}

\newcommand\gl{{\mathfrak{gl}}}

\newcommand\GL{{\mathrm{GL}}}

\def\Aut{{\mathrm {Aut}}}

\def\C{{\mathcal C}}
\def\D{{\mathcal D}}

\def\L{{\mathcal L}}

\def\O{{\mathcal O}}

\def\U{{\mathcal U}}

\def\Z{{\mathcal Z}}

\def\sdim{{\mathrm {sdim}}}
\def\ad{{\mathrm {ad}}}

\def\Ker{{\mathrm {Ker}}}

\def\Hom{{\mathrm {Hom}}}

\def\cHom{{\mathcal{H}om}}

\def\End{{\mathrm {End}}}
\def\cEnd{{\mathcal{E}nd}}

\def\Der{{\mathrm {Der}}}
\def\InnDer{{\mathrm {InnDer}}}

\def\ind{{\mathrm {ind}}}

\def\Cliff{{\mathrm {Cliff}}}

\def\a{{\mathfrak a}}
\def\g{{\mathfrak g}}
\def\h{{\mathfrak h}}
\def\k{{\mathfrak k}}
\def\n{{\mathfrak n}}

\def\st{{\mathfrak {st}}}

\def\p{{\mathfrak{p}}}

\def\z{{\mathfrak{z}}}

\def\place{{-}}

\begin{document}

\maketitle

\begin{abstract}
In this article we prove that there exists a Dixmier map for nilpotent super Lie algebras. 
In other words, if we denote by $\mathrm{Prim}(\U(\g))$ the set of (graded) primitive ideals of the enveloping algebra $\U(\g)$ of a nilpotent Lie superalgebra $\g$ 
and $\mathcal{A}d_{0}$ the adjoint group of $\g_{0}$, 
we prove that the usual Dixmier map for nilpotent Lie algebras can be naturally extended to the context of nilpotent super Lie algebras, \textit{i.e.} there exists a bijective map
\[     I : \g_{0}^{*}/\mathcal{A}d_{0} \rightarrow \mathrm{Prim}(\U(\g))     \] 
defined by sending the equivalence class $[\lambda]$ of a functional $\lambda$ to a primitive ideal $I(\lambda)$ of $\U(\g)$, 
and which coincides with the Dixmier map in the case of nilpotent Lie algebras. 
Moreover, the construction of the previous map is explicit, and more or less parallel to the one for Lie algebras, a major difference with a previous approach (\textit{cf.}~\cite{Let92}). 
One key fact in the construction is the existence of polarizations for super Lie algebras, generalizing the concept defined for Lie algebras. 
As a corollary of the previous description, we obtain the isomorphism $\U(\g)/I(\lambda) \simeq \Cliff_{q}(k) \otimes A_{p}(k)$, where $(p,q) = (\dim(\g_{0}/\g_{0}^{\lambda})/2,\dim(\g_{1}/\g_{1}^{\lambda}))$, we get a direct construction of the maximal ideals of the underlying algebra of $\U(\g)$ 
and also some properties of the stabilizers of the primitive ideals of $\U(\g)$. 
\end{abstract}



\section*{Introduction}

The aim of this article is to extend the Kirillov orbit method \textit{\`a la Dixmier} for nilpotent Lie algebras to the context of nilpotent super Lie algebras. 
More precisely, we shall prove the following results. 
Let $\g$ be a nilpotent super Lie algebra over an algebraically closed field of characteristic $0$. 
First, for every linear functional $\lambda \in \g_{0}^{*}$ there exists a so called polarization $\h$ of $\g$ at $\lambda$ (see Subsection \ref{subsec:superpol}) 
such that the induced module $\ind(\lambda|_{\h},\g)$ is simple and the kernel of its structure morphism 
is a (graded) primitive ideal of the universal enveloping algebra $\U(\g)$ (see Theorem \ref{teo:teo1}). 
Moreover, the previously constructed ideal does not depend on the polarization (see Theorem \ref{teo:teo2}), and it will be denoted $I(\lambda)$. 
Conversely, for every (graded) primitive ideal $I$ of the universal enveloping algebra $\U(\g)$, there exists a linear functional $\lambda \in \g_{0}^{*}$ 
such that $I = I(\lambda)$ (see Theorem \ref{teo:teo3}) and we further have that $\U(\g)/I \simeq \Cliff_{q}(k) \otimes A_{p}(k)$, 
where $(p,q) = (\dim(\g_{0}/\g_{0}^{\lambda})/2,\dim(\g_{1}/\g_{1}^{\lambda}))$ and $\g^{\lambda}= (\g^{\lambda}_{0},\g^{\lambda}_{1})$ is the kernel of the superantisymmetric bilinear form determined by $\lambda$ on $\g$ (see Proposition \ref{prop:mio}). 
A similar version of this last result was proved for any field of characteristic zero by A. Bell and I. Musson in \cite{BM90}, 
but without any determination of the indices $(p,q)$. 
Finally, we have that $I(\lambda) = I(\lambda')$ if and only if $\lambda$ and $\lambda'$ are in the same coadjoint orbit on $\g_{0}^{*}$ under the action of the adjoint group $\mathcal{A}d_{0}$ of $\g_{0}$ (see Proposition \ref{prop:6.2.3}). 
Summarizing, if we denote by $\mathrm{Prim}(\U(\g))$ the set of primitive ideals of $\U(\g)$, these results can be equivalently restated as saying that the map
\[     I : \g_{0}^{*}/\mathcal{A}d_{0} \rightarrow \mathrm{Prim}(\U(\g))     \] 
given by sending the equivalence class $[\lambda]$ of a functional $\lambda$ to $I(\lambda)$ is well-defined and bijective. 
As a consequence, we also derive an explicit description of the maximal ideals of the underlying algebra of $\U(\g)$. 
If $\g$ is just a nilpotent Lie algebra, the previous results are exactly the statements of the Dixmier map, 
which were gradually proved (together with generalizations to the solvable and other cases) by N. Conze, J. Dixmier, M. Duflo and M. Vergne to say a few names 
(\textit{cf.}~\cite{Dix96}, Ch. 6 and the references therein, especially the Supplementary remarks at \S 6.6). 

One of our main motivations to consider this extension is to study representations arising from (noncommutative) supersymmetric gauge field theory in physics, 
that could be found proceeding in an analogous manner to the one used for Yang-Mills theory in the Ph.D. thesis of the author, 
in which the Kirillov orbit method for nilpotent Lie algebras was used extensively (\textit{cf.}~\cite{HS10}, where these results were published). 
More precisely, I studied the representation theory of the so-called Yang-Mills algebras, defined by A. Connes and M. Dubois-Violette in \cite{CDV02}, and established a connection between them and the Weyl algebras, which are related to the gauge theory of the noncommutative flat space, using the Dixmier map. 
In this article we prove the necessary results of the Kirillov orbit method for nilpotent super Lie algebras to extend the previous point of view 
to a superized version, which is related to supersymmetric gauge field theory. 

We would like to make a few comments. 
At first glance, it could seem that what we have proved follows from the work \cite{Kac77} of V. Kac 
(\textit{cf.}~Thm. 7' on p. 82, where he claims similar results for completely solvable super Lie algebras). 
Nonetheless, as stated by A. Sergeev in \cite{Ser99}, Theorems 7 and 7' in \cite{Kac77} contain a mistake. 
In particular the proof of one key point of them, namely item (a) of the first of these theorems, does not hold, as explained by Sergeev in Sec. 5 of the aforementioned article. 
We would like to remark however that items (b), (c) and (d) of Theorem 7 still hold for nilpotent super Lie algebras, as it can be deduced from the results in \cite{Ser99}. 
On the other hand, we mention that E. Letzter has proved in \cite{Let92} that there is a bijection from the set of (graded) primitive ideals of the enveloping 
algebra of a completely solvable super Lie algebra to the set of primitive ideals of the enveloping algebra of the Lie algebra of even elements 
of the given super Lie algebra. 
This implies the existence of a bijection $I$ as before, being the inverse of the map studied by Letzter. 
However, his construction is not very explicit and in fact it does not provide a description of the ideals of the enveloping algebra of the super Lie algebra. 
We would like to remark that our manner of proceeding is completely different: even though we made use of his results as a short cut in the proof 
of Theorem \ref{teo:teo3} and Proposition \ref{prop:6.2.3} (for which a direct proof following the nonsuper case is possible), 
our aim is a more explicit description of the maximal ideals of $\U(\g)$ which relies on the structure of the super Lie algebra $\g$, 
following more or less the pattern of the nonsuper case. 
This last way of proceeding has allowed us, for instance, to study the quotients of the enveloping algebra of a nilpotent super Lie algebra by its primitive ideals, 
a description of the maximal ideals of the underlying algebra of the enveloping algebra of a nilpotent super Lie algebra 
and other properties related to our motivations (\textit{cf.}~Subsection \ref{subsec:sevcons}). 
Finally, we would also like to say it also generalizes some results proved by S. Mukherjee in \cite{Mu04} (\textit{cf.}~Thm. 11.1 of the mentioned paper). 

A lot of what we state at the beginning will be the superized version of properties already known for Lie algebras 
or associative algebras. 
Since many of the proofs of these result are more or less the same as for the nonsuper case, we shall only state the superized versions of them and put a standard reference 
whose proof can also be applied in the super world with at most minor changes, that will be shortly explained. 
Such changes could include the use of homogeneous elements instead of general ones, and the usually consequent use of super commutators instead of commutators, 
the appearance of inessential (and obvious) signs (\textit{cf.}~Lemma \ref{lema:4.5.1} for a typical case) and the replacement of the use of previous results for algebras by their superized versions for super algebras. 
This should not make the reader believe that in the graded world the standard theory of rings follows verbatim, as one can notice by reading \cite{NvO04} or \cite{CM84}. 
In the particular case of enveloping algebras of super Lie algebras, there were also several obstacles, as the construction of polarizations 
(because the super version of \cite{Dix96}, Prop. 1.12.10, does not yield polarizations of solvable super Lie algebras), the impossibility to exponentiate odd derivations in order to obtain automorphisms, etc. 
There are thus a collection of proofs that were not so clear at first sight, because they rely on other results inside a lengthy chain of generalizations, 
and others that were not clear at all, at least to the author. 
We shall include these. 

The article is organized as follows. 
In the first section we recall some generalities on super (associative) algebras, their ideals, representations and localizations. 
In order to fix notation, we also provide a short reminder on super vector spaces. 
In Section \ref{sec:lie}, we recall general facts on super Lie algebras and their representations. 
We also remind the usual definition of the universal enveloping algebra, and discuss some of its properties. 
In particular, we shall recall some basic facts of the ideals of the enveloping algebra of super Lie algebras 
and their representations that will be used in the construction of the Dixmier map. 
As already explained, in these two sections we shall provide a list of results that will be used in the sequel and that will be the corresponding superized version of well-known 
facts for plain algebras and Lie algebras. 

In Section \ref{sec:pol} we provide a definition of polarization for super Lie algebras, which resembles more to the definition for Lie algebras 
than to the one implicit in \cite{Kac77} for super Lie algebras. 
In order to do that, we first recall some definitions on bilinear forms on super vector spaces. 
Moreover, we remind the basic facts on polarizations of Lie algebras and some results proved by Sergeev in \cite{Ser99}. 
At the end, we give the main result of the section, namely the fact that solvable super Lie algebras have polarizations and derive some consequences. 
In the final section we state and prove the main results of this article, described at the beginning, and we derive several consequences. 

\section*{Acknowledgements}

I would like to thank J. Alev for having introduced me to the Kirillov orbit method and to the Dixmier map, and for extremely illuminating conversations. 
On the other hand, while achieving the major results of this article, M. Duflo has made me notice that he had previously worked on the same problem 
and obtained similar results, which were not published. 
I would like to thank him for his encouragement to write this article, for many interesting and enlightening conversations, and especially, for several comments about polarizations, 
and for pointing out reference \cite{BBB07}. 

\section{Generalities on super algebras}

From now on, we choose $k$ to be an algebraically closed field of characteristic $0$. 
All unadorned tensor products $\otimes$ are over $k$.

In the first subsection we recall the basic definitions of super vector spaces. 
Since the terminology sometimes varies in the literature, this section is also useful to fix the notation and some simple results we shall use in the sequel. 
We more or less follow the conventions of \cite{DM99}, which we suggest as a reference. 

In Subsection \ref{subsec:supalg} we will provide the definitions and some basic results on super algebras and their representations, 
stating several results concerning the two-sided ideals of a super algebra, that shall be crucial to us. 
Finally, in the last subsection we will remind some useful basic facts on localization of super algebras. 

\subsection{Super vector spaces}

We recall that a \emph{super vector space} over $k$ is a $k$-vector space $V$ provided with a $\ZZ/2\ZZ$-grading of the form $V = V_{0} \oplus V_{1}$. 
An element $v \in V$ is called \emph{homogeneous} if $v \in V_{i}$ for some $i \in \ZZ/2\ZZ$, and more precisely, 
the elements of $V_{0}$ are called \emph{even} and the ones which belong to $V_{1}$ are called \emph{odd}. 
For a nonzero homogeneous vector $v \in V_{i}$ ($i \in \ZZ/2\ZZ$), we write $|v| = i$ and call it the \emph{degree} or \emph{parity} of $v$. 
When we speak about the parity of an element, we will always assume that it is homogeneous. 

Given two super vector spaces $V$ and $W$, a \emph{morphism (of super vector spaces)} $f : V \rightarrow W$ is a $k$-linear map between the underlying vector spaces that preserves the grading. 
The vector space of morphisms from $V$ to $W$ is denoted by $\Hom(V,W)$.
It is easy to see that the collection of super vector spaces provided with the previous morphisms is a $k$-linear category, denoted by $\mathrm{sVect}_{k}$. 
A subobject of an object of this category will be called a \emph{sub super vector space} or more simply a \emph{subspace}.
We may define the \emph{super dimension} (called \emph{dimension} in \cite{DM99}) $\sdim(V)$ of a super vector space $V$ as the pair $(\dim(V_{0}),\dim(V_{1}))$. 

The category $\mathrm{sVect}_{k}$ is provided of a functor $\Pi$, called \emph{parity}, which satisfies that $\Pi(V)_{i} = V_{1-i}$, for $i \in \ZZ/2\ZZ$, and if $f : V \rightarrow W$, then $\Pi(f) = f$. 
Moreover, it is easy to see that this category is monoidal, when considering the tensor product given by $(V \otimes W)_{i} = \oplus_{j \in \ZZ/2\ZZ} V_{j} \otimes V_{i-j}$, for $i \in \ZZ/2\ZZ$, and the unit given by the super vector space of super dimension $(1,0)$, which we shall denote by $k$ (instead of $\underline{1}$ in \cite{DM99}). 
We may also consider the \emph{internal hom} in the category which is the super vector space $\cHom(V,W)$ (instead of $\underline{\Hom}$ in \cite{DM99}) 
such that $\cHom(V,W)_{0} = \Hom(V,W)$ 
and $\cHom(V,W)_{1} = \Hom(V,\Pi(W))$. 
It is clear that there is an adjunction between the tensor product and the internal hom of the form 
\[     \Hom(V \otimes W, U) \simeq \Hom(V,\cHom(W,U)).     \]
In fact, the previous isomorphism is just the degree zero part of a natural isomorphism of super vector spaces 
\[     \cHom(V \otimes W, U) \simeq \cHom(V,\cHom(W,U)).     \]
If we consider the flip $V \otimes W \rightarrow W \otimes V$ given by $v \otimes w \rightarrow (-1)^{|v||w|} w \otimes v$, 
one sees that $\mathrm{sVect}_{k}$ is in fact a braided monoidal category. 

\subsection{Super algebras}
\label{subsec:supalg}

A \emph{super associative and unitary algebra} is a vector space $A$ provided with morphisms of super vector spaces $\mu : A \otimes A \rightarrow A$, called \emph{product}, and 
an element $1_{A} \in A_{0}$ such that $\mu (\mu (a \otimes b) \otimes c)  =  \mu (a \otimes \mu(b \otimes c))$ for all $a, b, c \in A$ 
and $\mu (a \otimes 1_{A}) = a = \mu (1_{A} \otimes a)$, for all $a \in A$. 
As usual, we denote the product by a dot or simply by juxtaposition: $\mu (a \otimes b) = a \cdot b = a b$. 
By simplicity, super algebra will always denote a super associative and unitary algebra. 
A \emph{morphism} $\phi : A \rightarrow B$ of super algebras $A$ and $B$ is a morphism of the underlying super vector spaces $\phi : A \rightarrow B$ such that 
$\phi(a a')= \phi(a)\phi(a')$, for all $a, a' \in A$, and $\phi(1_{A}) = 1_{B}$. 
The \emph{tensor product} of super algebras is canonically defined following the Koszul's sign rule. 

Every super algebra is provided with an \emph{isomorphism} $\Sigma$ of order two, defined as $\Sigma(a_{0}+a_{1}) = a_{0} - a_{1}$, 
where $a_{i}$ is a homogeneous element of degree $i$, for $i \in \ZZ/2\ZZ$. 
For a super algebra $A$, we denote by $\O(A)$ its underlying algebra. 
From the previous comments, it is clear that a super algebra can be equivalently defined as an algebra provided with an isomorphism of algebras of order two. 

\begin{ex}
\label{ex:alg}
If $V$ is a super vector space, the internal endomorphism space $\cEnd(V)$ is a super algebra with the product given by composition. 
\end{ex}

A \emph{left module} of a super algebra $A$ is a super vector space $V$ provided with a morphism of super algebras $\rho : A \rightarrow \cEnd(V)$. 
Given two left $A$-modules $V$ and $W$, a \emph{morphism} from $V$ to $W$ is a morphism between the underlying super vector spaces $f : V \rightarrow W$ such that 
$f(a v) = a f(v)$, for all $a \in A$ and $v \in V$. 
There are similar definitions for right $A$-modules. 

The following examples of super algebras will be of great importance to us. 
\begin{ex}
\label{ex:weylcliff}
\begin{itemize}
\item[(i)] Given $n \in \NN_{0}$, let $A_{n}(k)$ denote the super algebra over $k$ given by 
\[     k\cl{q_{1},\dots,q_{n},p_{1},\dots,p_{n}}/\cl{\{[q_{i},p_{j}]-\delta_{i j}1, 1 \leq i,j \leq n\}},     \] 
where $q_{i}$ and $p_{i}$ are homogeneous of degree $0$, for $i= 1, \dots, n$ 
(we remark that $A_{0}(k) = k$). 
It is thus concentrated in degree zero, so a plain algebra, and it is called the $n$-th \emph{Weyl algebra}. 

\item[(ii)] Let $V$ be a vector space over $k$ with a nondegenerate symmetric bilinear form $\cl{\hskip 0.6mm,}$. 
Define $\Cliff(V,\cl{\hskip 0.6mm,})$ the super algebra given by $TV/\cl{\{v \otimes w + w \otimes v - \cl{v,w} 1 \hskip 0.6mm : \hskip 0.6mm v, w \in V \}}$, 
where the elements of $V$ are all of degree $1$. 
Then, it is a super algebra, called the \emph{Clifford algebra} of $(V,\cl{\hskip 0.6mm,})$. 
It can be proved that it only depends on the dimension $n$ of $V$, so it will be also 
denoted by $\Cliff_{n}(k)$. 
Moreover, it is easy to see that $\Cliff_{1}(k) \simeq k[\epsilon]/(\epsilon^{2}-1)$, where $|\epsilon| = 1$, and that $\Cliff_{2}(k) \simeq M_{2}(k)$, 
where the even and odd parts of the matrix algebra $M_{2}(k)$ are
\[     \begin{pmatrix} k & 0 \\ 0 & k \end{pmatrix} \hskip 0.5cm \text{and} \hskip 0.5cm \begin{pmatrix} 0 & k \\ k & 0 \end{pmatrix},     \]
respectively. 
Furthermore, we have the so called \emph{Bott periodicity}, \textit{i.e.} $\Cliff_{n+2}(k) \simeq \Cliff_{n}(k) \otimes M_{2}(k)$ (see \cite{Ka78}, Ch. III, Subsec. 3.23). 
We set $\Cliff_{0}(k) = k$. 
\end{itemize}
\end{ex}

The following proposition is known for algebras (see \cite{Dix63}). 
The proof for super algebras is more or less similar, but we include it for completeness. 
\begin{prop}
\label{prop:simplealg}
Suppose that $k$ is uncountable. 
Let $A$ be a super algebra and let $V$ be a simple module over $A$, which we assume to have a countable homogeneous basis over $k$. 
Then, every $A$-linear endomorphism of $V$ is given by a multiplication by a scalar in $k$. 
\end{prop}
\noindent\textbf{Proof.}
By Schur's Lemma in \cite{Ra98}, p. 591, we see that $\cEnd_{A}(V)$ is a super field (\textit{i.e.} every nonzero homogeneous element is invertible), 
so $\End_{A}(V)$ is a field. 
Following an idea of Dixmier, let us suppose that there exist an isomorphism $\phi \in \End_{A}(V)$ such that it is not the same as the multiplication by a scalar in $k$. 
This is equivalent to the fact that $\phi$ is not algebraic over $k$, for $k$ is algebraically closed. 
Hence, $\phi$ is trascendental over $k$ and $\End_{A}(V)$ contains a copy of the field $k(\phi)$, which has an uncountable basis, because $k$ is not countable. 
However, since $V$ is simple, any nonzero homogeneous element is a generator of the $A$-module $V$, say $v$, 
so every endomorphism of $V$ is completely determined by its value at $v$, which is a linear combination of the countable homogeneous basis of $V$. 
So the dimension of $\End_{A}(V)$ over $k$ is at most countable. 
This is a contradiction, so $\phi$ must be given by the multiplication by a scalar in $k$. 
The proposition is thus proved. 
\qed

\begin{rem}
\label{rem:simplealg}
We remark the fact that $\End_{A}(V) \simeq k$ immediately implies that $\cEnd_{A}(V)$ is $k$ or $k[\epsilon]/(\epsilon^{2}-1)$, with $|\epsilon| = 1$ 
(
see \cite{Va04}, Section 6.2, p. 215). 
This can be proved as follows. 
We first show that $\cEnd_{A}(V)_{1}$ has super dimension less than or equal to $1$. 
Let us suppose that it is not zero. 
We will prove that is one, \textit{i.e.} that two homogeneous isomorphisms in $\cEnd_{A}(V)$ degree $1$ are linearly dependent. 
Given $\phi, \psi \in \cEnd_{A}(V)$ two homogeneous isomorphisms of degree $1$, there exists a nonzero $c_{\phi,\psi} \in k$ such that $\phi \circ \psi = c_{\phi,\psi} 1_{V}$. 
From that we see that $c_{\phi,\phi} 1_{V} \circ \psi = \phi \circ \phi \circ \psi = \phi \circ c_{\phi,\psi} 1_{V}$. 
So, $\phi$ is just a scalar multiple in $k$ of $\psi$, and hence $\dim(\cEnd_{A}(V)_{1}) = 1$.   
Moreover, since $\phi^{2} = c_{\phi,\phi} 1_{V}$, with $c_{\phi,\phi} \neq 0$, we may define $\epsilon = \phi/\sqrt{c_{\phi,\phi}}$. 
\end{rem}

A \emph{subalgebra} of a super algebra $A$ is a subspace of the underlying super vector space of $A$ such that it is closed under the product of $A$. 
A \emph{left (resp. right, two-sided) ideal} of $A$ is subspace $I$ of the super vector space underlying $A$ such that $a x \in I$ (resp. $x a \in I$, $a x a' \in I$, for all $a' \in A$ and) for all $a \in A$ and $x \in I$. 
A two-sided ideal will be usually called ideal. 
For clarity, we remark that in this article, the term left (resp. right, two-sided) ideal of a super algebra will always denote 
a so-called left (resp. right, two-sided) \emph{graded} or \emph{super ideal}, 
which are sometimes used in the literature. 
Note however that we do distinguish between the ideals of a super algebra $A$ and the ideals of the underlying algebra $\O(A)$ of $A$. 

Given two homogeneous elements $a, b \in A$, the \emph{super commutator} $[a,b]$ of $a$ and $b$ is defined as $a b - (-1)^{|a||b|} b a$. 
We recall that a homogeneous element $z \in A$ is called \emph{supercentral} if $[z,a] = 0$, for all homogeneous elements $a \in A$. 
The \emph{super center} of $A$ is the super vector space expanded by the supercentral elements of $A$. 
A homogeneous $k$-linear map $d$ in $\cEnd(A)$ is called a \emph{derivation} if it satisfies the super Leibniz identity, \textit{i.e.} 
\[     d(a b) = d(a) b + (-1)^{|a||d|} a d(b).     \]

A super algebra $A$ is called \emph{left (resp. right) noetherian} if any left (resp. right) ideal has a finite set of homogeneous generators. 
Equivalently, $A$ is left (resp. right) noetherian if it satisfies the ascending chain condition on left (resp. right) ideals. 
From now on, noetherian will always denote left noetherian, unless we say the contrary. 
It is obvious to see that if $A$ is noetherian as an algebra, then it is noetherian as a super algebra. 

If $R$ is an algebra provided with an isomorphism $\sigma$ and a $\sigma$-derivation $\delta$ 
(\textit{i.e.} $\delta(r r') = \delta(r) r' + \sigma(r) \delta(r')$, for all $r, r' \in R$), 
$R[t,\sigma,\delta]$ will denote the \textit{Ore extension} of $R$, which is the (unique) algebra with underlying vector space $\oplus_{i \in \NN_{0}} R.t^{i}$, such that the product 
extends the left action of $R$ on the direct sum, $t^{i} t^{j} = t^{i+j}$ and $t r = \sigma(r) t + \delta(r)$, for all $r \in R$. 
Note that, if $R$ is the underlying graded of a super algebra $A$ and $\sigma = \Sigma$, then a $\sigma$-derivation is just an odd derivation of $A$. 
More generally, if $R$ is the underlying algebra of a super algebra $A$, $\sigma$ is an isomorphism of $A$ and $\delta$ is a homogeneous element of the internal morphism space of the underlying super vector space of $A$ of degree $|\delta|$, 
which is a $\sigma$-derivation of $\O(A)$, then $\O(A)[t,\sigma,\delta]$ has also the structure of a super algebra where $|t|=|\delta|$, 
which we shall denote by $A[t,\sigma,\delta]$. 

We shall now recall some properties of two-sided ideals of super algebras. 

An ideal $I$ of a super algebra $A$ is called \emph{maximal} if $I \neq A$ and it is maximal in the set of all ideals of $A$ different from $A$ with respect to inclusion. 
It is called \emph{primitive} if it is the annihilator of a simple left $A$-module. 
The \emph{(Jacobson) radical} $J(A)$ of $A$ is the intersection of all primitive ideals, or equivalently, the intersection of all maximal ideals of $A$ 
(\textit{cf.}~\cite{CM84}, Sec. 4, p. 250). 
Since $\mathrm{char}(k) \neq 2$, we have that the Jacobson radical of the super algebra $A$ coincides with the Jacobson radical of the underlying algebra of $A$ (\textit{cf.}~\cite{CM84}, Thm. 4.4, (3)).  

Moreover, $I$ is called \emph{prime} if $I \neq A$ and if whenever $J K \subseteq I$, for $J, K$ ideals of $A$, then $J \subseteq I$ or $K \subseteq I$. 
Equivalently, $I$ is prime if $I \neq A$ and for $a, b \in A$ homogeneous elements not in $I$, we have that $a A b \nsubseteq I$. 
The super algebra $A$ is \emph{integral} if $A \neq 0$ and the product of two nonzero homogeneous elements is nonzero. 
An ideal $I$ of the super algebra $A$ is \emph{completely prime} if $A/I$ is integral. 
The ideal $I$ is called \emph{semiprime} if $I \neq A$ and if in $A/I$ every two-sided nilpotent ideal is null. 
It is obvious that the intersection of an arbitrary collection of semiprime ideals is semiprime. 

It is trivial to see that a completely prime ideal is prime, and that a prime ideal is semiprime. 
The standard arguments show that a maximal ideal is primitive and that a primitive ideal is prime (\textit{cf.}~\cite{Dix96}, 3.1.6). 

There is a strong relation between the concept of prime or maximal ideal for a super algebra and the same notion for the underlying algebra. 
\begin{lema}[\cite{CM84}, Lemma 5.1 and Thm. 6.3, \textit{cf.}~\cite{BM90}, Lemma 1.2]
\label{lema:primax}
Let $A$ be a super algebra and $I$ an ideal of $A$. 
The following are equivalent:
\begin{itemize}

\item[(i)] $I$ is a prime (resp. maximal) ideal of the super algebra $A$.

\item[(ii)] $I = P \cap \Sigma(P)$, for some prime (resp. maximal) ideal $P$ of the underlying algebra of $A$.

\item[(iii)] $I$ is a semiprime ideal of the underlying algebra of $A$ and its minimal prime ideals form an orbit under $\Sigma$ 
(resp., and are maximal ideals of the underlying algebra of $A$). 
\end{itemize}
\end{lema}

We recall that a super algebra $A$ over $k$ is said to be \emph{central simple} if its super center is $k$ and it has no nontrivial two-sided ideals. 
We remark that we do not require $A$ to be semisimple, as in  \cite{Va04}, Section 6.2. 
As examples of central simple super algebras we have $A_{n}(k)$ (see \cite{FD93}, Part III, Exercise 26) and $\Cliff_{n}(k)$, for $n \in \NN_{0}$ (see \cite{Va04}, Section 6.2, p. 215). 
We refer to \cite{Lam80}, Ch. 4, \S 2, or \cite{Va04}, Section 6.2, for a more detailed study on (finite dimensional) central simple super algebras.
The following analogous result to the Azumaya-Nakayama's Theorem will be used in the sequel.
\begin{lema}
\label{lema:4.5.1}
Let $A$  be a central simple super algebra, $B$ a super algebra, $\mathcal{I}$ the set of ideals of $B$ and $\mathcal{I}'$ the set of ideals of $A \otimes B$. 
Then, the map from $\mathcal{I}$ to $\mathcal{I}'$ given by $I \mapsto A \otimes I$ is a bijection. 
Also, the super center of the tensor product is given by $\Z(A \otimes B) = \Z(B)$. 
Moreover, $I$ is a maximal (resp. prime) ideal of $B$ if and only if $A \otimes I$ is a maximal (resp. prime) ideal of $A \otimes B$. 
\end{lema}
\noindent\textbf{Proof.}
The proof of the first two statements is analogous to the one given in \cite{FD93}, Thm. 3.5 and Lemma 3.7, but 
taking into account that all elements must be homogeneous and one should use super commutators instead of commutators 
(\textit{cf.}~\cite{Lam80}, proof of Thm. 2.3).
This immediately implies the assertion concerning maximal ideals. 
The proof of the statement for prime ideals is the same as the one given for \cite{Dix96}, Lemma 4.5.1. 
\qed

As a direct corollary of the previous result we have (\textit{cf.}~\cite{Lam80}, Thm. 2.3):
\begin{coro}
\label{coro:fd}
Let $A$ and $B$ be two central simple super algebras. 
Then the tensor product $A \otimes B$ is also a central simple super algebra. 
\end{coro}

\subsection{Localization of super algebras} 

In this subsection we shall recall some facts on localization of super algebras. 
Even though some of these results may be stated in more general terms, we restrict ourselves to the cases we need. 
We refer to \cite{NvO04} for a more comprehensive exposition. 
Furthermore, some of the results we will state here are the obvious generalizations (with the standard proofs) of those that can be found for instance in 
\cite{Dix96}, Ch. 3, \S 6, for the case of algebras. 

If $S$ is a subset of homogeneous elements of $A$, it is said to satisfy the \emph{graded left (resp. right) Ore conditions} 
(or to \emph{allow of an arithmetic of left (resp. right) fractions}) if 
\begin{itemize}
\item[(i)] $1 \in S$, $0 \notin S$ and $S$ is multiplicative closed,
\item[(ii)]  If $a \in A$ is a homogeneous element and $s \in S$ are such that $a s = 0$ (resp. $s a =0$), then there exists $s' \in S$ such that $s' a = 0$ (resp. $a s' = 0$).
\item[(iii)] For $s \in S$ and $a \in A$ homogeneous elements, there exist $t \in S$  and $b \in A$ (resp. $t' \in S$ and $b' \in A$) such that $t a = b s$ (resp. $a t' = s b'$). 
\end{itemize}
We remark that, for a set $S$ of homogeneous elements, 
the graded left (resp. right) Ore conditions are equivalent to the usual left (resp. right) Ore conditions (see \cite{NvO04}, Lemma 8.1.1). 
Moreover, the left (resp. right) version of the Ore condition $(ii)$ is always satisfied if $A$ is left (resp. right) noetherian as an algebra. 
If both the left and right graded Ore conditions are satisfied, the super algebra $A$ together with its subset $S$ of homogeneous elements is said to 
\emph{allow of an arithmetic of fractions}. 
If the algebra $A$ allows of an arithmetic of fractions, the left and right localization rings $S^{-1}A$ and $AS^{-1}$ can be defined in the obvious way, 
they are super algebras and in fact coincide (see \cite{NvO04}, Prop. 8.1.2). 
We denote any of them by $A_{S}$. 
If $z$ is a homogeneous element and $S_{z}=\{z^{n} : n \in \NN_{0}\}$, we will usually write $A_{z}$ instead of $A_{S_{z}}$.
From now on we shall restrict to the case that $S$ does not contain zero divisors and that all of its elements are of degree $0$. 
Then, we have the following result:
\begin{lema}
Let $A$ be a super algebra and $S$ a set of elements of degree $0$ allowing of an arithmetic of fractions, and let $I$ be a two-sided ideal of $A$ satisfying that, 
if $s a \in I$, for $s \in S$ and $a \in A$, then $a \in I$. 
Let us consider $IS^{-1}$ and $S^{-1}I$ the subspaces of the super vector space underlying $A_{S}$ expanded by the homogeneous elements $a s^{-1}$ and $s^{-1} a$ respectively, 
where $a \in I$ is homogeneous and $s \in S$. 
Then, $S^{-1}I \subseteq IS^{-1}$ and $IS^{-1}$ is in fact a two-sided ideal of $A_{S}$. 
\end{lema}
\noindent\textbf{Proof.}
The standard proof given in \cite{Dix96}, Lemma 3.6.14 works in this case as well, taking into account that all elements there should be homogeneous. 
\qed

\begin{prop}
\label{prop:3.6.15}
Let $A$ be a super algebra and $S$ a set of elements of degree $0$ allowing of an arithmetic of fractions. 
Define $\mathcal{I}_{S}$ the set of two-sided ideals of $A_{S}$ and $\mathcal{I}$ the set of two-sided ideals of $A$ satisfying the following property: 
either $a s \in I$ or $s a \in I$, for $s \in S$ and $a \in A$ homogeneous, implies that $a \in I$. 
Then, if $I \in \mathcal{I}$, it holds that $S^{-1}I = IS^{-1}$, which we simply denote by $I_{S}$, and the maps from $\mathcal{I} \rightarrow \mathcal{I}_{S}$ given by 
$I \mapsto I_{S}$ and $\mathcal{I}_{S} \rightarrow \mathcal{I}$ by $I' \mapsto I' \cap A$ are mutually inverse. 
Furthermore, if $I \in \mathcal{I}$ is prime, so is $I_{S}$. 
\end{prop}
\noindent\textbf{Proof.}
The standard proof given in \cite{Dix96}, Prop. 3.6.15 also applies in this case, taking into account that all elements there should be homogeneous, and replacing the use of \cite{Dix96}, Lemma 3.6.14 by the previous lemma. 
\qed

Finally, we have the following simple result.
\begin{prop}
\label{prop:3.6.18}
Let $A$ be a super algebra and $S$ a set of elements of degree $0$ allowing of an arithmetic of fractions. 
Then any derivation $d : A \rightarrow A$ can be extended to a unique derivation $d_{S} : A_{S} \rightarrow A_{S}$. 
More precisely, if $s^{-1} x = y r^{-1}$, for $s, r \in S$ and homogeneous elements $x, y \in A$, then 
\[     d_{S}(s^{-1} x) = - s^{-1} d(s) s^{-1} x + s^{-1} d(x) = d(y) r^{-1} - (-1)^{|y|} y r^{-1} d(r) r^{-1}.     \]
\end{prop}
\noindent\textbf{Proof.}
The proof in \cite{Dix96}, Prop. 3.6.18 also applies in this case, taking into account that all elements there should be homogeneous and the appearance 
of harmless signs due to the Koszul's sign rule. 
\qed

\section{Generalities on super Lie algebras}
\label{sec:lie}

In the first subsection we shall provide the basic definitions and results on super Lie algebras, most of those can be found in \cite{Kac77} or \cite{Sch79}, which we suggest as a reference. 
We will also recall the standard relation between super Lie algebras and super algebras given by the universal enveloping algebra. 
In the next subsection we shall focus on representations of super Lie algebras, stating some results on induced modules that will be used all throughout the paper, 
and a basic fact on ideals of enveloping algebras of nilpotent super Lie algebras. 
Finally, in the next subsection we provide a useful criterion for a representation of an enveloping algebra of a super Lie algebra to be simple, 
analogous to the well-known one for plain Lie algebras, which will allow us to prove the primitivity of the ideals to be considered in Section \ref{sec:dixmap}. 

\subsection{Basic facts on super Lie algebras}

A \emph{super Lie algebra} over the field $k$ is a super vector space $\g  =\g_{0} \oplus \g_{1}$ provided with a morphism 
\[     [\hskip 0.6mm,] : \g \otimes \g \rightarrow \g,     \]
called the \emph{Lie bracket}, such that 
\begin{itemize}
\item the bracket is \emph{superskewsymmetric}; \textit{i.e.} $[x,y] = - (-1)^{|x||y|} [y,x]$, for all nonzero homogeneous $x, y \in \g$, 

\item the bracket satisfies the super Jacobi identity; \textit{i.e.} $[x,[y,z]] =[[x,y],z] + (-1)^{|x||y|} [y,[x,z]]$, for all nonzero homogeneous $x, y, z \in \g$. 
\end{itemize}
Instead of the most common denomination ``Lie superalgebra'', which appears in \cite{Kac77}, we prefer to use the more systematic terminology in \cite{DM99}.
From now on, even though it is not necessary in many definitions, we will suppose that the underlying vector space of the super Lie algebra is finite dimensional. 

A \emph{morphism} $\phi : \g \rightarrow \g'$ between two super Lie algebras $\g$ and $\g'$ is a morphism of the underlying super vector spaces 
$\g \rightarrow \g'$ such that $\phi([x,y])=[\phi(x),\phi(y)]$, for all homogeneous $x, y \in \g$. 
We thus have the category of super Lie algebras. 

The following proposition is direct.
\begin{prop}[\cite{Va04}, Section 3.1, p. 89]
A super vector space $\g$ is a super Lie algebra if and only if $\g_{0}$ is a Lie algebra, $\g_{1}$ is a module over $\g_{0}$, 
there exists a $\g_{0}$-equivariant linear map $S^{2}\g_{1} \rightarrow \g_{0}$, and it holds that 
\[       [x,[x,x]] = 0, \hskip 0.6cm \forall  \hskip 0.6mm x \in \g_{1}.       \]
\end{prop}

\begin{ex}
\label{ex:lie}
\begin{itemize}
\item[(i)] Given a super algebra $A$, we may regard it as a super Lie algebra with the bracket given by the \emph{super commutator} $[a,b] = a b - (-1)^{|a||b|} b a$. 
We denote this structure by $\mathrm{sLie}(A)$. 

\item[(ii)] Given a super vector space $V$ of super dimension $(n,m)$ we may consider the super vector space $\cEnd(V)$. 
It has the structure of a super algebra (explained in Example \ref{ex:alg}), and 
thus of a super Lie algebra provided with the Lie bracket given by the super commutator. 
It is denoted by $\gl(V)$ or $\gl(n|m)$. 
\end{itemize}
\end{ex}

Given a super Lie algebra $\g$, we can consider a super algebra associated to it, called the \emph{universal enveloping algebra} $\U(\g)$, 
which is defined as the quotient of the tensor algebra $T\g$ by the ideal generated by $\{ x \otimes y - (-1)^{|x||y|} y \otimes x - [x,y] \}$ for all homogeneous $x, y \in \g$. 
The $\ZZ/2\ZZ$-grading on $\U(\g)$ is induced from the $\ZZ/2\ZZ$-grading of $\g$. 
This super algebra satisfies the universal property $\Hom(\U(\g),A) \simeq \Hom(\g,\mathrm{sLie}(A))$, 
where the first morphism space is of super algebras and the second one is of super Lie algebras. 
It is provided with an increasing filtration of super vector spaces $\{F^{\bullet}\U(\g)\}_{\bullet \in \NN_{0}}$ coming from the filtration of the tensor algebra $T\g$ 
given by its usual grading. 
It is easy to prove that the underlying algebra of the super algebra $\U(\g)$, for $\g$ a finite dimensional super Lie algebra, is noetherian 
(see \cite{Be87}, Prop. 3.1, (i)), so \textit{a fortiori} the super algebra $\U(\g)$ is noetherian. 

As well as for enveloping algebras of Lie algebras, there is a PBW theorem for super Lie algebras. 
\begin{teo}[\cite{Ro65}, Thm. 2.1]
\label{teo:pbw}
Let $\{ x_{1}, \dots, x_{s} \}$ be an ordered basis of $\g$ consisting of homogeneous elements. 
Then the set of all products of the form 
\[     x_{1}^{p_{1}} \dots x_{s}^{p_{s}},   \]
where $x_{i}^{0} = 1$, $p_{i} \in \NN_{0}$ and $p_{i} \leq 1$ whenever $x_{i}$ is odd, is a basis for $\U(\g)$.
\end{teo}

We recall that an \emph{antiautomorphism} of a super algebra $A$ is an isomorphism $\phi$ of the underlying super vector space 
satisfying that $\phi(x y) = (-1)^{|x||y|} \phi(y) \phi(x)$, for all homogeneous elements $x, y \in A$, and $\phi(1)=1$. 
The enveloping algebra $\U(\g)$ is provided with an antiautomorphism $\alpha$, called \emph{principal}, such that $\alpha(x) = - x$, for $x \in \g$. 
In fact, it is easily proved that  
\[     \alpha(x_{1} \dots x_{n}) = (-1)^{n+ \sum_{i<j} |x_{i}||x_{j}|} x_{n} \dots x_{1},     \]
where $x_{1}, \dots, x_{n}$ are homogeneous elements of the super Lie algebra $\g$, and $n \in \NN$. 

A \emph{sub super Lie algebra} of a super Lie algebra $\g$ is a super vector space $\h \subseteq \g$ closed under the bracket operation, \textit{i.e.} if $x, y \in \h$, then $[x,y] \in \h$. 
Analogously, a \emph{super Lie ideal} of a super Lie algebra $\g$ is a super vector space $\k \subseteq \g$ that satisfies that, for all $x \in \g$ and $y \in \k$, $[x,y] \in \k$. 
Equivalently, we could have given the previous two definitions just in terms of homogeneous elements. 
Since we shall often work in the ``super'' context, we will usually use the shorter terms subalgebra and ideal, unless we need to make the distinction. 

Given a super vector space $V \subseteq \g$ (resp. a set of homogeneous elements $S \subseteq \g$), the \emph{super centralizer} of $V$ (resp. $S$) is the super vector space 
$\C(V)$ (resp. $\C(S)$) expanded by the homogeneous elements $x \in \g$ such that $[x,y]=0$, for all homogeneous $y \in V$ (resp. $y \in S$). 
It is easily seen to be a subalgebra of $\g$. 
The super centralizer of $\g$ is called the \emph{super center} of the super Lie algebra $\g$, and denoted by $\Z(\g)$. 
It is an ideal of $\g$. 

On the other hand, given two super vector spaces $V, W \subseteq \g$ (resp. two sets of homogeneous elements $S, T \subseteq \g$), the \emph{super commutator} $[V,W]$ (resp. $[S,T]$) 
of $V$ and $W$ (resp. of $S$ and $T$) is the super vector space expanded by $[x,y]$, for all the homogeneous elements $x \in V$ and $y \in W$ (resp. $x \in S$ and $y \in T$). 
If in the previous definition $V$ and $W$ are ideals, the super commutator is also an ideal. 
In particular, $[\g,\g]$ is called the \emph{derived algebra} of $\g$. 

\begin{ex}
\label{ex:der}
Given a super Lie algebra $\g$, the super vector space $\Der(\g)$ expanded by the homogeneous maps $d \in \cHom(\g,\g)$ satisfying that 
\[     d([x,y]) = [d(x),y] + (-1)^{|x||d|} [x,d(y)]     \]
is called the space of \emph{derivations} of $\g$. 
It is a super Lie algebra with the bracket provided by the super commutator. 
It is clear that the image $\InnDer(\g)$ of the morphism of super Lie algebras $\ad : \g \rightarrow \cHom(\g,\g)$ is an ideal of $\Der(\g)$, 
called the space of \emph{inner derivations} of $\g$. 
Note that $\ad$ induces an obvious identification $\InnDer(\g)_{0} = \InnDer(\g_{0})$. 
\end{ex}


As in the nongraded situation, we may consider the \emph{lower central series} of $\g$ to be the decreasing sequence of ideals defined recursively by 
$\C^{1}(\g)=\g$ and $\C^{i}(\g) = [\g,\C^{i-1}(\g)]$ for $i \geq 2$. 
Furthermore, the \emph{derived series} of $\g$ is the decreasing sequence of ideals defined recursively by $\D^{0}(\g)=\g$ and $\D^{i}(\g) = [\D^{i-1}(\g),\D^{i-1}(\g)]$ for $i \in \NN$. 
It is easy to see that $\D^{i}(\g) \subseteq \C^{i+1}(\g)$, for all $i \in \NN_{0}$. 
A super Lie algebra $\g$ is called \emph{solvable} if the there exists $i \in \NN_{0}$ such that $\D^{i}(\g) = 0$. 
Analogously, $\g$ is said to be \emph{nilpotent} if the there exists $i \in \NN_{0}$ such that $\C^{i}(\g) = 0$. 
It is clear that a nilpotent super Lie algebra is solvable. 

The following result indicates that the solvability of a super Lie algebra only relies on its even part.
\begin{prop}[\cite{Kac77}, Prop. 1.3.3, or \cite{Ser99}, Cor. 2.3]
\label{prop:sol}
A super Lie algebra $\g$ is solvable if and only if the Lie algebra $\g_{0}$ is solvable.
\end{prop}

There is also a version of Engel's theorem for super Lie algebras and it is proved in exactly the same way as for (plain) Lie algebras. 
\begin{prop}[\cite{Sch79}, Ch. III, \S 2, 1., Prop. 1]
\label{prop:engel}
Let $\g$ be a subalgebra of $\gl(n|m)$ such that, for every homogeneous $x \in \g$, $\mathrm{ad}(x)$ is a nilpotent operator. 
Then there is a homogeneous vector $v$ in the super vector space of super dimension $(n,m)$ such that $x(v)=0$, for all $x \in \g$.
\end{prop}

As a corollary of the previous result we have: 
\begin{coro}[\cite{Sch79}, Ch. III, \S 2, 1., Cor. 1]
\label{coro:nil}
A super Lie algebra $\g$ is nilpotent if and only if for every homogeneous $x \in \g$, $\mathrm{ad}(x)$ is a nilpotent operator. 
As a consequence, a super Lie algebra $\g$ is nilpotent if and only if $\g_{0}$ is a nilpotent Lie algebra and the action of $\g_{0}$ on $\g_{1}$ is by nilpotent operators. 
\end{coro}

Concerning the enveloping algebra of a nilpotent Lie super algebra, we now state a basic result that we will used throughout the article. 
\begin{prop}[\cite{Let89}, Prop. 3.3] 
Let $I$ be an ideal of the enveloping algebra $\U(\g)$ of a nilpotent super Lie algebra $\g$, distinct from $\U(\g)$. 
Then, $I$ is primitive if and only if it is maximal.
\end{prop}

\subsection{Representations of super Lie algebras}

A \emph{left representation} of a super Lie algebra $\g$ (or a \emph{left $\g$-representation}) is a super vector space $V$ provided with a morphism of super Lie algebras 
$\rho : \g \rightarrow \gl(V)$. 
Equivalently, a left representation of $\g$ is a super vector space $V$ provided with a morphism of super vector spaces 
\[     \rho' : \g \otimes V \rightarrow V     \]
such that, for all homogeneous $x, y \in \g$, 
\[     \rho'(x \otimes \rho'(y \otimes v)) - (-1)^{|x||y|} \rho'(y \otimes \rho'(x \otimes v)) = \rho'([x,y] \otimes v).     \]
It is clear that $\rho'(x \otimes v) = \rho(x)(v)$. 
We shall usually denote the action by a dot or even by juxtaposition, \textit{i.e.} $\rho(x)(v) = x \cdot v = x v$. 

Given a left $\g$-representation $V$ with structure morphism $\rho$, the parity changed representation $\Pi V$ is defined as follows. 
The underlying super vector space is just the parity functor $\Pi$ applied to the underlying super vector space of the $\g$-representation $V$. 
However, the action satisfies the identity $x.v = (-1)^{|x|} \rho(x)(v)$, for homogeneous $x \in \g$ and $v \in \Pi V$, and where the left member stands for the action of 
$x$ on $v \in \Pi V$, but on the right member we are considering the action of $x$ on $V$. 

Given two left representations $V$ and $W$ of $\g$, a \emph{morphism} $f : V \rightarrow W$ is a map of the underlying super vector spaces such that 
$f (x v) = x f(v)$. 
We denote the space of such morphisms by $\Hom_{\g}(V,W)$. 
We will also consider the super vector space of morphisms $\cHom_{\g}(V,W)$ given by $\Hom_{\g}(V,W)_{0}=\Hom_{\g}(V,W)$ and $\Hom_{\g}(V,W)_{1}=\Hom_{\g}(V,\Pi W)$. 
There are similar definitions for right representations. 


It is trivial to see that the category of left representations of $\g$ is equivalent to the category of left modules over the super algebra $\U(\g)$, since 
$\Hom(\g,\gl(V)) \simeq \Hom(\U(\g),\cEnd(V))$, where the first morphism space if of super Lie algebras and the second one is of super algebras. 
From now on, we will deal only with left representations and modules, and just call them representations and modules, respectively. 
Moreover, for a representation $V$ of $\g$, the corresponding morphisms $\g \rightarrow \gl(V)$ and $\U(\g) \rightarrow \cEnd(V)$ are called the \emph{structure morphisms} of $V$. 

\begin{ex}
The \emph{adjoint representation} of $\g$ in itself given by $x.y = \ad(x)(y) = [x,y]$ can be extended by derivations to a representation in $\U(\g)$, 
which is called the \emph{adjoint representation} of $\g$ in $\U(\g)$. 
We shall denote this representation by $\U(\g)^{\mathrm{ad}}$ and the structure morphism by $\ad$. 
More generally, if $\k$ is an ideal of $\g$, then it is a subrepresentation of the adjoint representation of $\g$ in itself, and it can also be extended by derivations 
to a representation of $\g$ in $\U(\k)$, which is also called the \emph{adjoint representation} of $\g$ in $\U(\k)$.
\end{ex}

The following lemma is the super version of \cite{Dix96}, Lemma 2.2.22 and will be needed later. 
The proof is similar to the nonsuper case, but we provide it because of the signs, which come from the use of the Koszul's sign rule. 
\begin{lema}[\textit{cf.}~\cite{Dix96}, Lemma 2.2.22]
\label{lema:auto}
Let $\k$ be an ideal of a super Lie algebra $\g$ and define $\delta = \ad \circ \alpha$, 
where $\alpha$ denotes the principal antiautomorphism of $\U(\g)$ and $\ad : \U(\g) \rightarrow \cEnd(\U(\k))$ 
is the structure morphism of the adjoint representation of $\g$ in $\U(\k)$. 
For $p \geq 0$, consider $y_{1}, \dots, y_{p} \in \g$ and $z \in \U(\k)$ homogeneous elements and $n_{1}, \dots, n_{p} \in \NN$, then 
\begin{equation}
\label{eq:auto}
      z y_{1}^{n_{1}} \dots y_{p}^{n_{p}} 
      = \sum_{0 \leq m_{i} \leq n_{i}} (-1)^{\epsilon_{m_{1},\dots,m_{p}}^{n_{1},\dots,n_{p}}(z)} \binom{n_{1}}{m_{1}} \dots \binom{n_{p}}{m_{p}} 
                                             y_{1}^{m_{1}} \dots y_{p}^{m_{p}} \delta(y_{1}^{n_{1}-m_{1}} \dots y_{p}^{n_{p}-m_{p}})(z),
\end{equation}
where 
\[     \epsilon_{m_{1},\dots,m_{p}}^{n_{1},\dots,n_{p}}(z) = |z|\sum_{i=1}^{p}n_{i}|y_{i}| + \sum_{i<j}(n_{i}-m_{i})|y_{i}|m_{j}|y_{j}|.     \]
\end{lema}
\noindent\textbf{Proof.}
We proceed by induction on $p$. 
The previous identity clearly holds for $p=0$. 
Let us assume that it is true for $p-1$, and we shall prove it for $p$. 
This implies that it holds for $p$ and $n_{1} = 0$. 
By induction again, we suppose that \eqref{eq:auto} holds for $p$ and $n_{1}-1$ (and arbitrary $n_{2}, \dots, n_{p}$), and we will prove it for $n_{1}$. 
Since $z y_{1} = (-1)^{|z||y_{1}|} (y_{1} z + \delta(y_{1})(z))$, we obtain that 
\begin{small}
\begin{align*}
 z y_{1}^{n_{1}} \dots y_{p}^{n_{p}} =& (-1)^{|z||y_{1}|} (y_{1} z + \delta(y_{1})(z)) y_{1}^{n_{1}-1} \dots y_{p}^{n_{p}} 
      \\
      =& (-1)^{|y_{1}||z|} \hskip -4mm
      \underset{\text{\begin{tiny}$    \begin{matrix}
                                       0 \leq m_{1} \leq n_{1}-1 \\ 0 \leq m_{i} \leq n_{i} \\ 2 \leq i \leq p
                                       \end{matrix}$
                                       \end{tiny}}}
                                       {\sum} \hskip -5mm
         (-1)^{\epsilon_{m_{1},\dots,m_{p}}^{n_{1}-1,\dots,n_{p}}(z)} \binom{n_{1}-1}{m_{1}} \dots \binom{n_{p}}{m_{p}} 
                                             y_{1}^{m_{1}+1} \dots y_{p}^{m_{p}} \delta(y_{1}^{n_{1}-m_{1}-1} \dots y_{p}^{n_{p}-m_{p}})(z)
      \\
      &+ (-1)^{\eta} \hskip -4mm
      \underset{\text{\begin{tiny}$    \begin{matrix}
                                       0 \leq m_{1} \leq n_{1}-1 \\ 0 \leq m_{i} \leq n_{i} \\ 2 \leq i \leq p
                                       \end{matrix}$
                                       \end{tiny}}}
                                       {\sum} \hskip -5mm
         (-1)^{\epsilon_{m_{1},\dots,m_{p}}^{n_{1}-1,\dots,n_{p}}(z)} \binom{n_{1}-1}{m_{1}} \dots \binom{n_{p}}{m_{p}} 
                                             y_{1}^{m_{1}} \dots y_{p}^{m_{p}} \delta(y_{1}^{n_{1}-m_{1}-1} \dots y_{p}^{n_{p}-m_{p}})(\delta(y_{1})(z))
      \\
      =& \sum_{0 \leq m_{i} \leq n_{i}} (-1)^{\epsilon_{m_{1},\dots,m_{p}}^{n_{1},\dots,n_{p}}(z)} \binom{n_{1}}{m_{1}} \dots \binom{n_{p}}{m_{p}} 
                                             y_{1}^{m_{1}} \dots y_{p}^{m_{p}} \delta(y_{1}^{n_{1}-m_{1}} \dots y_{p}^{n_{p}-m_{p}})(z),
\end{align*}
\end{small}
with $\eta=|y_{1}||z|+|y_{1}|((n_{1}-1)|y_{1}| + \sum_{i=2}^{p}n_{i}|y_{i}|)$. 
We remark that we have used the inductive assumption in the third member, and the identities 
\[     \delta(y_{1}^{n_{1}-m_{1}} \dots y_{p}^{n_{p}-m_{p}}) = (-1)^{|y_{1}|((n_{1}-m_{1}-1) |y_{1}| + \sum_{i=2}^{p}(n_{i}-m_{i})|y_{i}|)} 
\delta(y_{1}^{n_{1}-m_{1}-1} \dots y_{p}^{n_{p}-m_{p}}) \delta(y_{1})     \] 
and
\[     \binom{n_{1}-1}{m_{1}-1} + \binom{n_{1}-1}{m_{1}} = \binom{n_{1}}{m_{1}}     \]
in the last member.
\qed

A \emph{subrepresentation} of a representation $V$ is a subspace $W$ of the super vector space $V$ such that $\rho(x)(w) \in W$, for all $w \in W$ and all $x \in \g$. 
It is clear that any nonzero representation $V$ has at least two different subrepresentations: $V$ and $0$, which are called \emph{trivial}. 
We shall say that a representation $V$ is \emph{irreducible} or \emph{simple} if its only subrepresentations are trivial. 
We remark that an irreducible representation of the super Lie algebra $\g$, or equivalently, of the super algebra $\U(\g)$, 
may have nontrivial subspaces of the underlying vector space of $V$ which are invariant under the action of $\g$. 

\begin{prop}[\textit{cf.}~\cite{Dix96}, Prop. 2.6.5]
\label{prop:simple}
Let $V$ be a $\g$-representation. 
Then the following are equivalent
\begin{itemize}
 \item[(i)] $V$ is simple.

 \item[(ii)] $V$ is simple and every $\g$-linear endomorphism of $V$ is given by the multiplication by a scalar in $k$.

 \item[(iii)] $V \neq 0$ and, for any set of homogeneous elements $x_{1}, \dots, x_{n}, y_{1}, \dots, y_{n} \in V$ satisfying that $x_{i} \in V_{b}$ and $y_{i} \in V_{a}$ 
for fixed $a,b \in \ZZ/2\ZZ$ and for all $i=1, \dots , n$, and $x_{1}, \dots, x_{n}$ linearly independent over $k$, there exists a homogeneous $z \in \U(\g)$ such that $z x_{i} = y_{i}$, for all $i=1, \dots , n$.
\end{itemize}
\end{prop}
\noindent\textbf{Proof.}
The implication $(iii) \Rightarrow (i)$ (and also $(ii) \Rightarrow (i)$) is trivial. 
On the other hand, the implication $(ii) \Rightarrow (iii)$ follows from the Density Theorem for graded rings (\textit{cf.}~\cite{Ra98}, Lemma 2, and \cite{ELS04}, Thm. 1.3). 
Finally, the implication $(i) \Rightarrow (ii)$ can be proved following the lines of the classical proof given by D. Quillen in \cite{Qui69}, Thm. 1. 
We first note that the Generic Flatness Lemma holds by trivial reasons in a little more general context: 
\begin{lema}
\label{lema:gfl}
Let $A$ be a finitely generated commutative $k$-algebra, which is also domain, and let $B$ be an $A$-algebra, \textit{i.e.} we assume that there is a morphism 
of algebras $A \rightarrow \mathcal{Z}(B)$. 
We suppose further that $\mathcal{Z}(B)$ is a finitely generated $A$-algebra, under the previous morphism, and that $B$ is a finitely generated $\mathcal{Z}(B)$-module. 
Hence, if $N$ is a finitely generated $B$-module, there exists $a \in A$ such that $A_{a} \otimes_{A} N$ is free over the localization $A_{a}$ 
of $A$ at the element $a$. 
\end{lema}
\noindent\textbf{Proof.}
As stated before, the proof follows from the Generic Flatness Lemma (see \cite{Gro63}, 60--61, Expos\'e IV, Lemme 6.7) applied to the $A$-algebra 
$\mathcal{Z}(B)$, because any finitely generated $B$-module $N$ is also finitely generated when considered as a $\mathcal{Z}(B)$-module due to the 
assumption on $B$. 
\qed

The result now follows using the same argument given in the proof of the first theorem in \cite{Qui69} (\textit{cf.} also \cite{Dix96}, Lemma 2.6.4), 
which we repeat just for convenience. 
It suffices to prove that every nonzero element $\theta \in \End_{\U(\g)}(V)$ is algebraic over $k$. 
Since $V$ is simple, $\theta$ is an isomorphism. 
Set $A = k[\theta]$, which is a finitely generated commutative domain, for $\theta$ is an (even) isomorphism of the simple module $V$. 
It is a super algebra with trivial odd homogeneous component. 
We suppose that $\theta$ is not algebraic over $k$, so $A$ is a polynomial algebra. 
So, $V$ has the structure of a module over the tensor product super algebra $A \otimes \U(\g)$. 

Consider now $\U(\g)$ provided with the canonical filtration $F^{\bullet}\U(\g)$, whose associated graded algebra is $S(\g_{0}) \otimes \Lambda \g_{1}$. 
Set $v \in V$ be a nonzero element, and define an exhaustive filtration $F^{\bullet}V$ on $V$ compatible with the one of the enveloping algebra given by $F^{\bullet}V = A F^{\bullet}\U(\g) v$. 
Taking into account that $\mathrm{gr}_{F^{\bullet}V}(V)$ is a finitely generated module over the algebra $B = A \otimes \mathrm{gr}_{F^{\bullet}\U(\g)}(\U(\g))$, which satisfies the hypotheses 
of Lemma \ref{lema:gfl}, there exists $a \in k[\theta]$ such that $\mathrm{gr}_{F^{\bullet}V}(V)_{a}$ is free over $k[\theta]_{a}$. 
The module $V_{a}$ over $A_{a} \otimes \U(\g)$ is also provided with a compatible exhaustive filtration of the form $F^{\bullet}V_{a} = (F^{\bullet}V)_{a}$. 
The exactness of localization tells us that $\mathrm{gr}_{F^{\bullet}V}(V)_{a} \simeq \mathrm{gr}_{F^{\bullet}V_{a}}(V_{a})$ is a free $A_{a}$-module. 
Since $A_{a}$ is principal, each component $F^{p}V_{a}/F^{p-1}V_{a}$ is free over $A_{a}$, because it is finitely generated without torsion, 
and the fact that the filtration is bounded below and exhaustive yields that $V_{a}$ is also a free $A_{a}$-module. 

Consider $a' \in A$ such that $a'$ does not divide any power of $a$, so the multiplication by $a'$ induces a nonsurjective morphism of $A$-modules on $A_{a}$. 
Hence, the endomorphism of $A_{a}$-modules of $V_{a} = V \otimes_{A} A_{a}$ given by $v \otimes a'' \mapsto v \otimes a' a'' = a'(v) \otimes a''$ is not surjective. 
On the other hand, since $V$ is a module over the division ring $\Hom_{\U(\g)}(V,V)$, the map $v \mapsto a'(v)$ is an isomorphism, which gives a contradiction, and proves the claim. 
\qed

Given a subalgebra $\h$ of a super Lie algebra $\g$ and a representation $W$ of $\h$, the \emph{representation of $\g$ induced by $W$}, denoted by $\ind(W,\g)$, is given by 
$\U(\g) \otimes_{\U(h)} W$ with the left action given by the regular action of $\g$ on $\U(\g)$. 
It is clear that, given any $\g$-representation $V$, there is a canonical isomorphism 
\begin{equation}
\label{eq:iso}
     \Hom_{\h}(W,V) \simeq \Hom_{\g}(\ind(W,\g),V),
\end{equation}
where in the first morphism space we regard $V$ as an $\h$-representation coming from the inclusion $\h \subseteq \g$. 
The fact that $\U(\g)$ is free over $\U(\h)$ tells us that $W$ is simple if $\ind(W,\g)$ is simple. 

We recall that, given a representation $V$ of $\g$, the \emph{annihilator} of a subspace $W$ (resp. a subset $W$ of homogeneous elements) of the underlying 
super vector space of $V$ is the left ideal of $\U(\g)$ formed by the elements $x \in \U(\g)$ such that $x w = 0$, for all $w \in W$. 
\begin{prop}
\label{prop:5.1.7}
Let $\h$ be a subalgebra of a super Lie algebra $\g$ and let $W$ be a representation of $\h$ given by $\rho : \U(\h) \rightarrow \cEnd(W)$. 
Set $V=\ind(W,\g)$ the representation of $\g$ induced by $W$ and $J = \ker(\rho)$. 
Then,
\begin{itemize}
 \item[(i)] The annihilator of $W$ in $\U(\g)$ is the left ideal $\U(\g)J$. 
 
 \item[(ii)] The kernel of the induced structure morphism $\pi : \U(\g) \rightarrow \cEnd(V)$ is the largest two-sided ideal of $\U(\g)$ contained in $\U(\g)J$.
\end{itemize}
\end{prop}
\noindent\textbf{Proof.}
The proof given in \cite{Dix96}, Prop. 5.1.7, works word for word.
\qed

The following two results are easy but we state them just for clarity. 
\begin{lema}
\label{lema:5.1.8}
Let $\lambda \in \Hom(\g, k)$ be a functional of $\g$ such that $\lambda([\g,\g])=0$. 
It determines two one-dimensional representations of $\g$ of the form $F_{\lambda,\g,i} = k.v_{\lambda,\g,i}$ with $|v_{\lambda,\g,i}|=i$, 
for $i \in \ZZ/2 \ZZ$. 
The action is given by $x.v_{\lambda,\g,i} = \lambda(x) v_{\lambda,\g,i}$, for all $x \in \g$, and the structure morphism 
will also be donted by $\lambda$. 
Then, the kernel of the induced structure morphism $\U(\g) \rightarrow k$ coincides with the left (resp. right, two-sided) ideal of $\U(\g)$ 
generated by the set $\{ z - \lambda(z) \}$, for all homogeneous elements $z \in \g$. 
\end{lema}
\noindent\textbf{Proof.}
The proof given in \cite{Dix96}, Lemma 5.1.8, applies as well to this case.
\qed

\begin{prop}
\label{prop:5.1.9}
Let $\h$ be a subalgebra of a super Lie algebra $\g$ and let $W$ be a representation of $\h$ given by $\rho : \U(\h) \rightarrow \cEnd(W)$. 
Suppose further that $W$ has a (homogeneous) generator $w$ as an $\U(\h)$-module. 
Set $V=\ind(W,\g)$ the representation of $\g$ induced by $W$, and $L$ the annihilator of $w$ in $\U(\h)$. 
Then,
\begin{itemize}
 \item[(i)] The map from $\U(\g)$ to $V$ given by $u \mapsto u w$, for $u \in \U(\g)$ is surjective of kernel $\U(\g) L$, so it induces an isomorphism 
            of the form $\U(\g)/(\U(\g) L) \simeq W$. 
 
 \item[(ii)] If $W$ is one-dimensional, so the structure morphism can be written as $\rho : \U(\h) \rightarrow k$, then $\U(\g) L$ is the left ideal of $\U(\g)$ 
             generated by the set $\{ z - \lambda(z) \}$, for all homogeneous elements $z \in \h$.
\end{itemize}
\end{prop}
\noindent\textbf{Proof.}
The proof given in \cite{Dix96}, Prop. 5.1.8, works word for word, replacing the use of \cite{Dix96}, Lemma 5.1.9, by the previous lemma. 
\qed

The next result will be used in the sequel. 
\begin{prop}
\label{prop:5.1.13}
Let $\h$ be a subalgebra of a super Lie algebra $\g$, $\k$ an ideal of $\g$ contained in $\h$ and $W$ a representation of $\h$ given by $\rho : \U(\h) \rightarrow \cEnd(W)$ 
such that $\rho([\g,\k]) = 0$. 
Let $V = \ind(W,\g)$ be the representation of $\g$ induced by $W$ with structure morphism $\pi : \U(\g) \rightarrow \cEnd(V)$.
Then the $\k$-representation on $V$ given by $\pi|{\k}$ is a direct sum of copies of the $\k$-representation on $W$ given by $\rho|_{\k}$.
\end{prop}
\noindent\textbf{Proof.}
The proof follows the lines of that in \cite{Dix96}, Prop. 5.1.13, but with some changes. 
One first proves that, given $p \in \NN_{0}$ and homogeneous elements $y \in \k$ and $x_{1}, \dots, x_{p} \in \g$, we have that 
\[     y x_{1} \dots x_{p} \in (-1)^{|y|(|x_{1}|+\dots+|x_{p}|)} x_{1} \dots x_{p} y + \U(\g) [\g,\k].     \]
This is done by induction on $p$. 
It is obvious for $p = 0$. 
Let us suppose that it holds for $p-1$, and consider 
\begin{align*}
   y x_{1} \dots x_{p} &= (-1)^{|x_{1}||y|} x_{1} y x_{2} \dots x_{p} + [y,x_{1}] x_{2} \dots x_{p} 
   \\
   &\in (-1)^{|y|(|x_{1}|+\dots+|x_{p}|)} x_{1} \dots x_{p} y + x_{1} \U(\g) [\g,\k] 
   \\
   &+ (-1)^{(|y|+|x_{1}|)(|x_{2}|+\dots+|x_{p}|)} x_{2} \dots x_{p} [y,x_{1}] + \U(\g) [\g,\k]
   \\
   &\subseteq (-1)^{|y|(|x_{1}|+\dots+|x_{p}|)} x_{1} \dots x_{p} y + \U(\g) [\g,\k],     
\end{align*}
which establishes the claim. 

Now, if we consider homogeneous elements $y \in \k$, $u \in \U(\g)$ and $w \in W$, the previous result tells us that 
$y u = (-1)^{|u||y|} u y + \sum_{i \in I} u_{i} c_{i}$, for $u_{i} \in \U(\g)$ and $c_{i} \in [\g,\k]$. 
This yields that 
\[     \pi(y) (u \otimes_{\U(\h)} w) = (-1)^{|u||y|} u y \otimes_{\U(\h)} w + \sum_{i \in I} u_{i} \rho(c_{i}) \otimes_{\U(\h)} w 
         = (-1)^{|u||y|} u \otimes_{\U(\h)} \rho(y) w,     \]
for $\rho([\g,\k])=0$. 
If we consider a homogeneous basis $\{ z_{j} \}_{j \in J}$ of $\U(\g)$ over $\U(\h)$, we have that the $\k$-representation on $V$ given by $\pi|_{\k}$ is a 
direct sum of the $\k$-representations $k.z_{j} \otimes W$, where $\pi(y)(z_{j} \otimes w) = (-1)^{|z_{j}||y|} z_{j} \otimes \rho(y)w$. 
It is clear that $k.z_{j} \otimes W$ is isomorphic to $W$ if $z_{j}$ is even and to $W^{\Sigma}$ if $z_{j}$ is odd, where $W^{\Sigma}$ denotes the 
$\k$-representation on $W$ with structure morphism $\rho \circ \Sigma$. 
Finally, we note that the map of super vector spaces $W \rightarrow W^{\Sigma}$ given by $w_{0} + w_{1} \mapsto w_{0} - w_{1}$ is an isomorphism of 
$\k$-representations. 
The proposition is thus proved. 
\qed

We recall that a Lie algebra is called \emph{algebraic} if it is the Lie algebra of an algebraic group. 
Also, given $\g_{0}$ a finite dimensional Lie algebra, the set of automorphisms $\Aut(\g_{0})$ of the Lie algebra $\g_{0}$ is a linear algebraic group 
with (algebraic) Lie algebra $\Der(\g_{0})$ (see \cite{TY05}, Prop. 24.3.7). 
The algebraic action of $\Aut(\g_{0})$ on itself by conjugation induces a morphism of algebraic groups 
$\mathrm{Ad} : \Aut(\g_{0}) \rightarrow \GL(\Der(\g_{0}))$ (see \cite{TY05}, 23.5.2). 
Let $\mathfrak{ad}_{0}$ be the smallest algebraic subalgebra of the Lie algebra $\Der(\g_{0})$ satisfying that $\InnDer(\g_{0}) \subseteq \mathfrak{ad}_{0}$. 
The \emph{adjoint (algebraic) group} $\mathcal{A}d_{0}$ of $\g_{0}$ is the smallest algebraic subgroup of $\GL(\g_{0})$ (or $\Aut(\g_{0})$) whose Lie algebra contains $\InnDer(\g_{0})$, or equivalently, it is the irreducible algebraic subgroup of $\GL(\g_{0})$ (or $\Aut(\g_{0})$) with Lie algebra $\mathfrak{ad}_{0}$ (see \cite{TY05}, 24.8.1-2). 
The action of $\Aut(\g_{0})$ on $\Der(\g_{0})$ preserves the ideal $\InnDer(\g_{0})$, and thus preserves also the ideal $\mathfrak{ad}_{0}$. 
Hence, the action of $\Aut(\g_{0})$ on itself by conjugations preserves $\mathcal{A}d_{0}$. 
Note that $\mathcal{A}d_{0}$ acts naturally on $\g_{0}$. 
When $\InnDer(\g_{0}) = \mathfrak{ad}_{0}$, we will simply say that $\mathcal{A}d_{0}$ is the \emph{adjoint group} of $\g_{0}$. 
The previous identity is equivalent to say that $\g_{0}$ is algebraic, which is satisfied for every nilpotent Lie algebra.

The preceding paragraph can be extended to the case of algebraic super groups. 
Let now $\mathfrak{ad}'$ be the smallest algebraic subalgebra of the Lie algebra $\Der(\g)_{0}$ (or $\gl(\g)_{0}$) satisfying that $\InnDer(\g_{0}) \subseteq \mathfrak{ad}'$, and let 
$\mathcal{A}d'$ be the irreducible algebraic subgroup of $\Aut(\g)$ with Lie algebra $\mathfrak{ad}'$. 
The latter can be equivalently defined as the smallest algebraic subgroup of $\GL(\g_{0}) \times \GL(\g_{1})$ such that its Lie algebra contains the Lie algebra $\InnDer(\g_{0})$ 
(see \cite{TY05}, 24.8.1-2).
This immediately implies that $\mathfrak{ad}' = \mathfrak{ad}_{0}$ and $\mathcal{A}d' = \mathcal{A}d_{0}$ (see \cite{TY05}, Prop. 28.4.5). 
Arguing as before, we see that the action of $\Aut(\g)$ on $\Der(\g)_{0}$ preserves the ideal $\InnDer(\g_{0})$, 
and thus preserves also the ideal $\mathfrak{ad}_{0}$. 
So, the action of $\Aut(\g)$ on itself by conjugations preserves $\mathcal{A}d_{0}$ and we also have that $\mathcal{A}d_{0}$ acts naturally on $\g$.  


Let us suppose that we further have an action of an algebraic group $H_{0}$ on a super Lie algebra $\g$ by automorphisms of super Lie algebras, 
\textit{i.e.} that there is a morphism of algebraic groups $H_{0} \rightarrow \Aut(\g)$. 
It is easy to see that this action induces an action of $H_{0}$ on $\U(\g)$, which preserves the filtration of the enveloping algebra. 
This can be applied, to the standard action of the adjoint algebraic group $\mathcal{A}d_{0}$ of $\g_{0}$ on $\g$ given in the previous paragraph 
to see that it induces an action on $\U(\g)$. 

The following result will be used in the sequel. 
\begin{prop}
\label{prop:2.4.17}
Let $\g$ be a super Lie algebra and $\mathcal{A}d_{0}$ the adjoint algebraic group of $\g_{0}$. 
For every element $a \in \Aut(\g)$, denote by $a_{\U}$ the induced automorphism of the enveloping algebra $\U(\g)$. 
For any ideal $I$ of the enveloping algebra $\U(\g)$ and $a \in \mathcal{A}d_{0}$, we have that $a_{\U}(I)=I$. 
\end{prop}
\noindent\textbf{Proof.} 
The proof given in \cite{Dix96}, Prop. 2.4.17, applies word for word. 
\qed


\subsection{On the simplicity of representations of super Lie algebras}

Let $\k$ be an ideal of a super Lie algebra $\g$ and let $U$ be a representation of $\k$ with structure morphism $\sigma : \U(\k) \rightarrow \cEnd(U)$. 
Following \cite{Dix96}, the \emph{stabilizer of $\sigma$ in $\g$} is the subspace of the super vector space underlying $\g$ expanded by the homogeneous elements $y \in \g$ 
satisfying that there exists a homogeneous endomorphism $s \in \cEnd(U)$ of the same degree as $y$ such that
\[     \sigma([y,x]) = [s,\sigma(x)]     \]
for all homogeneous elements $x \in \k$, and where we remark that $[s,\sigma(x)]$ is the super commutator in $\cEnd(U)$. 
It is denoted by $\st(\sigma,\g)$ or $\st(U,\g)$. 
Analogously, given an ideal $I$ of $\U(\k)$ we may define the \emph{stabilizer of $I$ in $\g$}, denoted by $\st(I,\g)$, 
as the super vector space expanded by the homogeneous elements $x \in \g$ such that $\ad(x)(I) \subseteq I$ 
(\textit{i.e.} such that $\ad(x)(z) \in I$ for all homogeneous elements $z \in I$).
It is clear that both $\st(\sigma,\g)$ and $\st(I,\g)$ are subalgebras of $\g$ containing $\k$. 
Moreover, it is easy to see that if $\k$ is an ideal of a super Lie algebra $\g$, $U$ a representation of $\k$ such that its structure morphism $\sigma : \U(\k) \rightarrow \cEnd(U)$ 
has kernel $I$, then $\st(\sigma,\g) \subseteq \st(I,\g)$ (\textit{cf.}~\cite{Dix96}, Prop. 5.3.3). 

\begin{lema}
\label{lema:5.3.5}
Let $\k$ be an ideal of a super Lie algebra $\g$, $U$ a simple representation of $\k$ with structure morphism $\sigma : \U(\k) \rightarrow \cEnd(U)$, 
$\h = \st(\sigma,\g)$, $\rho : \U(\h) \rightarrow \cEnd(W)$ a representation of $\h$ such that the $\k$-representation on $W$ given by $\rho|_{\k}$ is a direct sum 
of copies of the $\k$-representation on $U$ given by $\sigma$, and let $V$ be the induced representation $\ind(W,\g)$ with structure morphism $\pi$. 
Let $V_{n}$ be the super vector space expanded by the classes of $x \otimes w$, for homogeneous elements $x \in F^{n}\U(\g)$ and $w \in W$.
It is an exhaustive increasing filtration of the super vector space $V$. 
Given $n \in \NN$ and $t \in V_{n} \setminus \{ 0 \}$, there exists $z \in \U(\k)$ such that $z t \in V_{n-1} \setminus \{0 \}$. 
\end{lema}
\noindent\textbf{Proof.}
The proof follows the pattern for the nonsuper case given in \cite{Dix96}, Prop. 5.3.5, but since there are several differences we give it. 

Let $\{ x_{1}, \dots x_{m}\}$ be a homogeneous basis of a complement of $\h$ in $\g$ and write 
$t = \sum_{|\bar{n}| \leq p} \bar{x}^{\bar{n}} \otimes_{\U(\h)} w_{\bar{n}}$, where $\bar{x}^{\bar{n}} = x_{1}^{n_{1}} \dots x_{m}^{n_{m}}$ and $w_{\bar{n}} \in W$. 
If $w_{\bar{n}} = 0$ for all $\bar{n}$ such that $|\bar{n}| = \sum_{i=1}^{m} n_{i} = p$, there is nothing to prove. 
It suffices thus to prove the lemma for the case that there exists some $\bar{n}_{0}$ with $|n_{0}| = p$ such that $w_{\bar{n}_{0}}$ does not vanish. 
We may suppose that the homogeneous element $t$ further satisfies that $t = \sum_{|\bar{n}| = p} \bar{x}^{\bar{n}} \otimes_{\U(\h)} w_{\bar{n}}$, 
where all nonvanishing $w_{\bar{n}} \in W$ are homogeneous of the same degree. 
This can be proved as follows. 
First, since $\U(\k)$ preserves the filtration defined on $V$, we may ignore the terms indexed by $\bar{n}$ with $|\bar{n}| < p$. 
Second, if we write $t = t^{0} + t^{1}$, where $t^{i}$ is the sum of the terms $\bar{x}^{\bar{n}} \otimes_{\U(\h)} w_{\bar{n}}$ such that $|w_{\bar{n}}|=i$, 
for $i \in \ZZ/2 \ZZ$, then the fact that the action of $\U(\k)$ preserves the filtration defined on $V$ implies that we may proceed stepwise, as we wanted to prove. 

Since $W$ can be written as a direct sum $\oplus_{\lambda \in \Lambda} W_{\lambda}$, for $W_{\lambda}$ a $\k$-representation isomorphic to $U$, let 
$\zeta_{\lambda} : W \rightarrow U$ be the unique epimorphism of $\k$-modules with kernel $\oplus_{\lambda' \neq \lambda} W_{\lambda'}$. 
Choose a nonzero homogeneous element $u \in U$ of the same parity as $w_{\bar{n}_{0}}$. 
By Proposition \ref{prop:simple}, there exists an even element $z \in \U(\k)$ and elements $\xi_{\lambda,w_{\bar{n}}} \in k$ 
such that $z \zeta_{\lambda} (w_{\bar{n}}) = \xi_{\lambda,w_{\bar{n}}} u$, for all $\lambda$ and $\bar{n}$. 
Lemma \ref{lema:auto} tells us that $z t \equiv \sum_{|\bar{n}|=p} \bar{x}^{\bar{n}} \otimes_{\U(\h)} z w_{\bar{n}}$ (\textrm{mod}. $V_{p-1}$), 
so by changing $t$ by $z t$, we can further assume that $t$ satisfies that $\zeta_{\lambda}(w_{\bar{n}}) = \xi_{\lambda,\bar{n}} u$, 
for some $\xi_{\lambda,\bar{n}} \in k$, and that there exists $\lambda_{0}$ such that $\xi_{\lambda_{0},\bar{n}_{0}} \neq 0$. 

Choose $i_{0} \in \{ 1, \dots, m \}$ such that $\bar{n}_{0} = (n_{0,1},\dots,n_{0,m})$ satisfies that $n_{0,i_{0}} \neq 0$, 
and define $\bar{n}'_{0} = \bar{n}_{0} - e_{i_{0}}$, where $e_{i_{0}}$ is the vector of $\NN_{0}^{m}$ which has $1$ in the $i_{0}$-th place and zero elsewhere. 
Given a homogeneous element $z \in \U(\k)$, Lemma \ref{lema:auto} yields that
\begin{equation}
\label{eq:1}
     z t \equiv \sum_{|\bar{n}|=p} (-1)^{|z||\bar{x}^{\bar{n}}|} \bar{x}^{\bar{n}} \otimes_{\U(\h)} z w_{\bar{n}} 
       - \sum_{j=1}^{m} \sum_{|\bar{n}|=p} (-1)^{|z||\bar{x}^{\bar{n}}|+|x_{j}|\sum_{i>j}n_{i}|x_{i}|} n_{j} \bar{x}^{\bar{n}-e_{j}} \otimes_{\U(\h)} 
      [x_{j},z] w_{\bar{n}} 
       \text{(\textrm{mod}. $V_{p-2}$)}.     
\end{equation}
Suppose now that the statement of the lemma does not hold, \textit{i.e.} that $\U(\k) t \cap V_{p-1} = 0$. 
If $z u = 0$, we get that $z t \in V_{p-1}$ and by the assumption it must vanish. 
Hence, we conclude that the component coming from $\bar{x}^{\bar{n}'_{0}}$ in \eqref{eq:1} is given by  
\[     \sum_{j=1}^{m} (-1)^{|x_{j}|(|z|+\sum_{i>j}n'_{0,i}|x_{i}|)} (n'_{0,j}+1) [x_{j},z] w_{\bar{n}'_{0}+e_{j}} = 0.     \]
By applying $\zeta_{\lambda_{0}}$ to this equality we obtain 
\[     \sum_{j=1}^{m} (-1)^{|x_{j}|(|z|+\sum_{i>j}n'_{0,i}|x_{i}|)} (n'_{0,j}+1) [x_{j},z] \xi_{\lambda_{0},\bar{n}'_{0}+e_{j}} u = 0,     \]
which can be rewritten as $[y,z]u=0$, for 
$y = \sum_{j=1}^{m} (-1)^{|x_{j}|(|z|+\sum_{i>j}n'_{0,i}|x_{i}|)} (n'_{0,j}+1)  \xi_{\lambda_{0},\bar{n}'_{0}+e_{j}} x_{j} \in \g$. 
It is direct to see that, if $y = y_{0} + y_{1}$ is the decomposition of $y$ in homogeneous elements, 
the previous vanishing identity is equivalent to $[y_{i},z]u=0$, for $i \in \ZZ/2\ZZ$. 
Define $y'_{0} = y_{0}$, $y'_{1} = (-1)^{|z|} y_{1}$ and $y'=y'_{0}+y'_{1}$. 
Moreover, the homogeneous component $y''$ of $y'$ of degree $|x_{i_{0}}|$ does not belong to $\h$ since $\xi_{\lambda_{0},\bar{n}_{0}} \neq 0$ and it does not depend on $z$. 

Let us define $s \in \End(U)$ given by $s(z u) = [y'',z] u$, for $z \in \U(\k)$. 
We recall that the degree of $s$ coincides with the degree of $y''$. 
It is well-defined by the previous considerations. 
Given homogeneous elements $z \in \U(\k)$ and $x \in \k$, we obtain that 
\[     [s,x](z u) = s (x z u) - (-1)^{|s||x|} x s(z u) = [y'', x z] u - (-1)^{|y''||x|} x [y'',z] u = [y'' , x] z u,     \]
which means that $y'' \in \h$, that is a contradiction. 
The lemma is thus proved. 
\qed

\begin{teo}
\label{teo:5.3.6}
Let $\k$ be an ideal of a super Lie algebra $\g$, $U$ a simple representation of $\k$ with structure morphism $\sigma : \U(\k) \rightarrow \cEnd(U)$, 
$\h = \st(\sigma,\g)$, and $\rho : \U(\h) \rightarrow \cEnd(W)$ a representation of $\h$ such that the $\k$-representation on $W$ given by $\rho|_{\k}$ is a direct sum 
of copies of the $\k$-representation on $U$ given by $\sigma$. 
Then the induced representation $\ind(W,\g)$ is simple. 
\end{teo}
\noindent\textbf{Proof.}
The proof for the nonsuper case in \cite{Dix96}, Thm. 5.3.6, works word for word in this case, replacing the use of \cite{Dix96}, Lemma 5.3.5 by 
Lemma \ref{lema:5.3.5}.
\qed

\section{Polarizations}
\label{sec:pol}

The aim of this section is to prove that polarizations exist for solvable super Lie algebras. 
In order to do so, we first provide some easy results on bilinear forms on super vector spaces. 
Then, we recall the basic facts on polarizations of Lie algebras. 
Finally, we will recall some of the ideas of A. Sergeev used to study irreducible finite dimensional representations of solvable super Lie algebras. 
As a consequence, we shall derive that all solvable super Lie algebras have polarizations. 
We would like to remark that M. Duflo has proved this result using a different idea (\textit{cf}. \cite{BBB07}). 

\subsection{Bilinear forms on super vector spaces}
\label{subsec:bil}

Let $V$ be a super vector space provided with an \emph{even bilinear form} $\cl{\hskip 0.6mm,}$, \textit{i.e.} a morphism of super vector spaces 
$\cl{\hskip 0.6mm,}: V \otimes V \rightarrow k$. 
We remark that the homogeneity of the map $\cl{\hskip 0.6mm,}$ is equivalent to the fact that $\cl{v,w} = 0$, for all $v, w \in V$ of different parity. 
We suppose moreover that $\cl{\hskip 0.6mm,}$ is either \emph{superantisymmetric} or \emph{supersymmetric}, \textit{i.e.} $\cl{v,w} = - (-1)^{|v||w|} \cl{w,v}$, or $\cl{v,w} = (-1)^{|v||w|} \cl{w,v}$, for all $v, w \in V$ homogeneous, respectively. 
Furthermore, we see that an even superantisymmetric (resp. supersymmetric) bilinear form on $V$ is equivalent to give an antisymmetric (resp. symmetric) bilinear form on $V_{0}$ and a symmetric (resp. antisymmetric) bilinear form on  $V_{1}$. 
From now on, all bilinear forms will be even unless otherwise stated. 

Let $W$ be a subspace of $V$. 
It is easy to prove that 
\[     \sdim(W) + \sdim(W^{\perp_{\cl{\hskip 0.1mm,}}}) = \sdim(V) + \sdim(W \cap V^{\perp_{\cl{\hskip 0.1mm,}}}),     \]
where $W^{\perp_{\cl{\hskip 0.1mm,}}}$ denotes the subspace of $V$ perpendicular to $W$ with respect to the form $\cl{\hskip 0.6mm,}$: 
\[     W^{\perp_{\cl{\hskip 0.1mm,}}} = \{ v \in V : \cl{v,w} = 0 \}.     \]
We remark that $W^{\perp_{\cl{\hskip 0.1mm,}}} = W_{0}^{\perp_{\cl{\hskip 0.1mm,}|_{V_{0}}}} \oplus W_{1}^{\perp_{\cl{\hskip 0.1mm,}|_{V_{1}}}}$. 
From now on, unless it is necessary to specify the super vector space and its form $\cl{\hskip 0.6mm,}$, 
we denote the perpendicular space to a subspace $W$ only by $W^{\perp}$. 
The subspace $V^{\perp}$ of the super vector space $V$ provided with $\cl{\hskip 0.6mm,}$ is called the \emph{kernel} of the bilinear form. 

We recall that $W \subseteq V$ is called \emph{totally isotropic} if $W \subseteq W^{\perp}$. 
Moreover, a totally isotropic subspace $W \subseteq V$ is \emph{maximal totally isotropic} if it is maximal in the set of totally isotropic subspaces of the super vector space $V$ 
with respect to the inclusion.  
We note that this implies that $W \subseteq W^{\perp}$, but it does not necessarily yield that $W = W^{\perp}$. 
All the previous definitions specialize to the usual ones for a vector space $V$ with subspace $W$ if we consider $V$ as a super vector space with $V_{1} = 0$.

It is easy to see that a subspace $W$ of a super vector space $V$ provided with a superantisymmetric or supersymmetric bilinear form $\cl{\hskip 0.6mm,}$ is (maximal) totally isotropic if and only if each $W_{i}$ is a (maximal) totally isotropic subspace of the vector space $V_{i}$ provided with $\cl{\hskip 0.6mm,}|_{V_{i}}$, 
for $i \in \ZZ/2\ZZ$. 
This implies that we may thus restrict to the study of symmetric and antisymmetric forms on vector spaces.

If $V$ is a vector space provided with an antisymmetric bilinear form and $W$ is a totally isotropic subspace then the following conditions are equivalent 
(see \cite{Dix96}, 1.12.1):
\begin{itemize}
\item $W$ is maximal in the set of totally isotropic subspaces with respect to the inclusion,

\item $\dim(W) = (\dim(V) + \dim(V^{\perp}))/2$,

\item $W \supset W^{\perp}$,

\item $W = W^{\perp}$.
\end{itemize}
On the other hand, let us assume that $V$ is a vector space with a symmetric bilinear form.
Taking the quotient by $V^{\perp}$, we may then restrict to the situation where the form is nondegenerate. 
In this case, a totally isotropic subspace $W$ of $V$ is maximal totally isotropic if and only if $W = W^{\perp}$, for $\dim(V)$ even, and $\dim(W) = \dim(W^{\perp})-1$, 
for $\dim(V)$ odd. 
Hence, for a vector space $V$ with an antisymmetric or symmetric bilinear form, the dimensions of all maximal totally isotropic subspaces coincide. 
In consequence, the super dimensions of all maximal totally isotropic subspaces of a super vector space coincide. 

\subsection{Polarizations of Lie algebras}
\label{subsec:pollie}

Let us first state the standard results about polarizations of plain Lie algebras.

A subalgebra $\h_{0}$ of a Lie algebra $\g_{0}$ is said to be \emph{subordinate} to a functional $\lambda_{0} \in \g_{0}^{*}$ if 
$\lambda_{0}([\h_{0},\h_{0}]) = 0$ (\textit{cf.}~\cite{Dix96}, 1.12.7). 
Equivalently, $\h_{0}$ is a totally isotropic subspace of $\g_{0}$ provided with the alternating bilinear form $A_{\lambda_{0}}$ given by 
$v \otimes w \mapsto \lambda_{0}([v,w])$.  
Moreover, we say that $\h_{0}$ is a \emph{polarization} of $\g_{0}$ at $\lambda_{0}$ if it is a subalgebra of $\g_{0}$ and it is a maximal totally isotropic subspace of the vector space underlying $\g_{0}$ provided with $A_{\lambda_{0}}$ (\textit{cf.}~\cite{Dix96}, 1.12.8). 
By the previous subsection, if $\g_{0}^{\lambda_{0}}$ denotes the kernel of $A_{\lambda_{0}}$, to be a maximally totally isotropic subspace is the same as to be totally isotropic and of dimension $(\dim(\g_{0}) + \dim(\g_{0}^{\lambda_{0}}))/2$.

Proposition 1.12.10 in \cite{Dix96} implies that, given any linear functional $\lambda_{0}$ on any solvable Lie algebra $\g_{0}$, a polarization of $\g_{0}$ at $\lambda_{0}$ always exists (we remark that one requires the assumption that $k$ is algebraically closed). 

\begin{rem}
\label{rem:pol}
We point out the easy fact that if $\k_{0}$ is an ideal of a Lie algebra $\g_{0}$ on which a functional $\lambda_{0} \in \g_{0}^{*}$ vanishes, 
then $\k_{0}$ should be included in every polarization of $\g_{0}$ at $\lambda_{0}$. 
Indeed, if $\h_{0}$ is a polarization at $\lambda_{0}$, then $\h_{0} + \k_{0}$ is also a subordinate subalgebra of $\lambda_{0}$. 
The maximality of $\h_{0}$ implies that $\k_{0} \subseteq \h_{0}$. 
\end{rem}

\subsection{Polarizations of super Lie algebras in the sense of Sergeev}

We shall now recall some definitions and facts from the work of Sergeev. 

Let $\g$ be a solvable super Lie algebra. 
Define $L_{\g}$ to be the vector space of functionals given by the elements $\lambda \in \Hom(\g, k)$ such that $\lambda([\g_{0},\g_{0}]) = 0$. 
We remark that the condition $\lambda \in \Hom(\g, k)$ is equivalent to say that $\lambda : \g \rightarrow k$ is a $k$-linear map between the underlying vector spaces 
such that $\lambda(\g_{1}) = 0$.
Analogously to the case of Lie algebras, a functional as before determines a symmetric bilinear form $B_{\lambda} : \g_{1} \otimes \g_{1} \rightarrow k$ given by 
$B_{\lambda}(x,y) = \lambda([x,y])$. 
A \emph{polarization in the sense of Sergeev} of $\g$ at $\lambda \in L_{\g}$ is a subalgebra $\h$ of $\g$ such that $\h_{0} = \g_{0}$ 
and $\h_{1}$ is a maximal totally isotropic subspace for the symmetric bilinear form $B_{\lambda}$. 

\begin{rem}
\label{rem:resser}
Note that, if $\h$ is a subspace of the super vector space underlying the super Lie algebra $\g$ such that $\h_{0} = \g_{0}$, 
then $\h$ is a subalgebra of $\g$ if and only if $\h_{1}$ is a $\g_{0}$-submodule of $\g_{1}$. 
\end{rem}

\begin{lema}[\cite{Ser99}, Lemma 2.4]
\label{lema:ser0}
Let $\g_{0}$ be a solvable Lie algebra and $V$ be a finite dimensional $\g_{0}$-module provided with a $\g_{0}$-invariant symmetric bilinear form. 
Given $W$ a $\g_{0}$-submodule of $V$, which is totally isotropic with respect to the bilinear form, then there exists a $\g_{0}$-submodule of $V$
which is a maximal totally isotropic subspace containing $W$. 
\end{lema}

The previous result applied to the case $V = \g_{1}$ and $W = 0$ (and taking into account Remark \ref{rem:resser}) implies: 
\begin{lema}[\cite{Ser99}, Lemma 1.1]
\label{lema:ser}
Let $\g$ be a solvable super Lie algebra. 
Given any $\lambda \in L_{\g}$, there exists a polarization in the sense of Sergeev of $\g$ at $\lambda$.
\end{lema}

\subsection{Polarizations of super Lie algebras}
\label{subsec:superpol}

Let us define $\L_{\g}$ to be the vector space of functionals given by the elements $\lambda \in \Hom(\g, k)$. 
Notice that $L_{\g} \subseteq \L_{\g}$. 
Given $\lambda \in \L_{\g}$, it defines a superantisymmetric bilinear form $\cl{\hskip 0.6mm,}_{\lambda}$ on $\g$ by the formula $\cl{x,y} = \lambda([x,y])$, for $x,y \in \g$, 
so \textit{a fortiori} an antisymmetric bilinear form $A_{\lambda}$ on $\g_{0}$ and a symmetric bilinear form $B_{\lambda}$ on $\g_{1}$. 
Denote by $\g^{\lambda}$ the kernel of $\cl{\hskip 0.6mm, }_{\lambda}$, which is a subalgebra of $\g$. 

V. Kac in \cite{Kac77}, p. 83, has defined a subalgebra $\h$ of $\g$ to be \emph{subordinate} to $\lambda$ if $\lambda([\h,\h])=0$ and $\h \supset \g^{\lambda}$. 
We define a \emph{polarization} of $\g$ at $\lambda \in \L_{\g}$ to be a subordinate subalgebra $\h$ of $\g$ such that it is a maximal totally isotropic subspace of the super vector space $\g$ with respect to the bilinear form $\cl{\hskip 0.6mm,}_{\lambda}$. 
By the considerations given in Subsection \ref{subsec:bil}, we see that $\h_{0}$ should be a polarization of the Lie algebra $\g_{0}$ at $\lambda|_{\g_{0}}$ and the super dimension of all polarizations at $\lambda$ coincide. 

We recall from Lemma \ref{lema:5.1.8} that $\lambda$ defines two one-dimensional $\h$-representations $F_{\lambda,\h,i} = k.v_{\lambda,\h,i}$ 
with $|v_{\lambda,\h,i}|=i$, for $i \in \ZZ/2 \ZZ$, with action given by $x.v_{\lambda,\h,i} = \lambda(x) v_{\lambda,\h,i}$, for all $x \in \h$. 
The structure morphism of this representation is also denoted by $\lambda|_{\h}$. 

\begin{rem}
If the super Lie algebra $\g$ is just a Lie algebra, this definition obviously coincides with the classical one given on Subsection \ref{subsec:pollie}. 
On the other hand, if $\lambda \in L_{\g}$, a subalgebra $\h \subseteq \g$ is a polarization at $\lambda$ if and only if it is a polarization in the sense of Sergeev at $\lambda$. 
This tells us that the new definition of polarization is an extension of the previous ones. 
\end{rem}

\begin{prop}
\label{prop:pol}
Let $\g$ be a solvable super Lie algebra. 
Every functional $\lambda \in \L_{\g}$ has a polarization at $\lambda$. 
\end{prop}
\noindent\textbf{Proof.}
Let $V \subseteq \g_{1}$ be a $\g_{0}$-submodule such that it is maximal totally isotropic with respect to the symmetric bilinear form $B_{\lambda}$. 
Such a submodule exists due to Lemma \ref{lema:ser0}. 
Now, consider the super commutator $[V,V] \subseteq \g_{0}$. 
It is easy verified that $[V,V]$ is a Lie ideal of the Lie algebra $\g_{0}$. 
Since $\g_{0}$ is a solvable Lie algebra by Proposition \ref{prop:sol}, there should exist a polarization $\p$ of $\g_{0}$ at $\lambda|_{\g_{0}}$. 
By definition of $V$, we have that $\lambda([V,V]) = 0$, and so, by Remark \ref{rem:pol} the Lie ideal $[V,V]$ of $\g_{0}$ should be included in any polarization of $\g_{0}$ at $\lambda|_{\g_{0}}$, and in particular $[V,V] \subseteq \p$. 
Since $V$ is a $\g_{0}$-submodule of $\g_{1}$, it is \textit{a fortiori} also an $\p$-submodule. 
This implies that the subspace $\h$ of $\g$ defined as $\h_{0} = \p$ and $\h_{1} = V$ is a subalgebra of $\g$. 
By construction, $\h_{i}$ is a maximal totally isotropic subspace of $\g_{i}$ provided with $\cl{\hskip 0.6mm,}|_{\g_{i}}$, for $i \in \ZZ/2\ZZ$. 
Hence, $\h$ is a polarization of $\g$ at $\lambda$. 
\qed

\begin{rem}
\label{rem:superpol}
Note that Remark \ref{rem:pol} also extends to this situation: if $\k$ is an ideal of a super Lie algebra $\g$ on which a functional 
$\lambda \in \L_{\g}$ vanishes, then $\k$ is included in every polarization of $\g$ at $\lambda$. 
The proof given there extends to this case word for word. 
In fact, the previous proof does not need that $\lambda$ vanishes on $\k$, only that it vanishes on $[\g,\k]$. 
\end{rem}

The following is a result of M. Duflo. 
We reproduce his proof because it does not seem to appear elsewhere. 
\begin{lema}[\textit{cf}. \cite{BBB07}, Cor. 5.2]
Let $\g$ be a solvable super Lie algebra and $V$ be a finite dimensional $\g$-module provided with a $\g$-invariant even superantisymmetric or supersymmetric 
bilinear form. 
Given $W$ a $\g$-submodule of $V$, which is totally isotropic with respect to the bilinear form, there exists a $\g$-submodule of $V$
which is a maximal totally isotropic subspace and contains $W$. 
\end{lema}
\noindent\textbf{Proof.}
Without loss of generality, let us suppose that the bilinear form $B$ on $V$ is superantisymmetric. 
If $B$ is supersymmetric then, we may consider $\Pi V$ instead of $V$. 
Moreover, taking $V/W$ instead of $V$, we may assume that $W = 0$. 

Define the super vector space $\h = V_{0} \oplus V_{1} \oplus k.z$, where we regard the homogeneous elements of $V$ with the same degree as in $V$ and $z$ in even degree. 
It is easy to see that $\h$ is a super Lie algebra if we define $[v,v'] = B(v,v') z$, for all $v, v' \in V$, and we declare $z$ to be supercentral. 
Furthermore, the action of $\g$ on $V$ extends to an action on $\h$ by derivations, if $\g.z = 0$, 
so we may consider the super Lie algebra given by the semidirect product $\g \ltimes \h$, and since $z$ is even, 
the functional $\lambda$ given by $\lambda|_{\g \oplus V} = 0$ and $\lambda(z) = 1$ belongs to $\L_{\g \ltimes \h}$. 
Since the semidirect product of solvable super Lie algebras is also solvable, we see that $\g \ltimes \h$ is solvable.
Finally, it is clear that the polarizations of $\g \ltimes \h$ at $\lambda$ are in bijection (taking the intersection with $V$) with the $\g$-modules 
of the statement. 
\qed

\begin{rem}
\label{rem:polinv}
We may apply the lemma to the following situations: 
\begin{itemize}
 \item[(i)] Given an ideal $\k$ of a solvable super Lie algebra $\g$, and a functional $\lambda \in \L_{\k}$, 
 then the previous result (for the $\g$-invariant superantisymmetric even bilinear form $\cl{\hskip 0.6mm, }_{\lambda}$ on $\k$) implies that there exist a polarization of 
 $\k$ at $\lambda$ that is invariant under $\g$.
 
 \item[(ii)] More generally, given solvable super Lie algebras $\h$ and $\k$ such that $\h$ acts by derivations on $\k$, and 
  a functional $\lambda \in \L_{\k}$, consider the super Lie algebra given by the semidirect product $\g = \h \ltimes \k$, which is obviously solvable. 
  Since $\k$ seen inside of $\g$ is an ideal, the previous item implies that there exists a polarization of $\k$ at $\lambda$ invariant under the action of $\h$ 
  (\textit{cf.}~\cite{Dix96}, Prop. 1.12.10, (iii)). 
\end{itemize}
\end{rem}

\section{The Dixmier map for nilpotent super Lie algebras}
\label{sec:dixmap}

We are now in position to prove the main results stated at the beginning of the introduction. 

\subsection{The main theorems}
\label{subsec:maindixmap}

We first recall the following result.
\begin{lema}[\cite{BM90}, Lemmas 1.10 and 2.2]
\label{lema:bm}
Let $\g$ be a nilpotent super Lie algebra with super center $\z = k z \neq \g$, where $z$ is even. 
Then, there exist elements $x, y \in \g$ homogeneous of the same parity and an ideal $\k$ of codimension one in $\g$ such that 
\begin{itemize}
 \item[(i)] $[y,x] = z$,
 \item[(ii)] $\k$ is the super centralizer of $y$ in $\g$ and $y \in \Z(\g/\z)$,
 \item[(iii)] $\g = \k \oplus k x$.
\end{itemize}
Moreover, if the super center of $\g/\z$ consists only of odd elements (so any $x$ and $y$ as before should be odd), then either of the following holds 
\begin{itemize}
 \item[(1)] $\g = k z \oplus k y$, with $[y,y]=z$ (\textit{i.e.} $x=y$),
 \item[(2)] there exists $y$ such that $(i)$, $(ii)$ and $(iii)$ hold and $[y,y]=0$. 
\end{itemize}
\end{lema}

\begin{lema}
\label{lema:abe}
A nilpotent super Lie algebra $\g$ of super dimension $(1,1)$ is isomorphic to one of the following:
\begin{itemize}
\item[(i)] $\g$ is supercommutative,

\item[(ii)] $\g = k z \oplus k c$, with $|z|=0$, $|c|=1$, $z \in \Z(\g)$ and $[c,c]=z$.
\end{itemize}
\end{lema}
\noindent\textbf{Proof.}
Let us suppose that $\g = k z \oplus k c$, with $|z|=0$, $|c|=1$.
By Corollary \ref{coro:nil}, $z$ must be supercentral. 
Then, the possibilities $(i)$ and $(ii)$ are equivalent to $[c,c]=0$ or $[c,c] \neq 0$, and the lemma follows. 
\qed

The following lemma will be useful when dealing with polarizations in an inductive process. 
\begin{lema}
\label{lema:mio}
Let $\g$ be a super Lie algebra with super center $\z = k z \neq \g$, where $z$ is even, and let $x, y \in \g$ be homogeneous of the same parity and $\k$ be an ideal 
of codimension one in $\g$ satisfying the properties (i), (ii), (iii) and (2) stated in Lemma \ref{lema:bm}. 
Given $\lambda \in \L_{\g}$ such that $\lambda(z) = 1$, define $\lambda' = \lambda|_{\k} \in \L_{\k}$.  
Then, if $\h$ is a polarization of $\k$ at $\lambda'$, it is also a polarization of $\g$ at $\lambda$. 
\end{lema}
\noindent\textbf{Proof.}
It is obvious that $\h$ is subordinate to $\lambda$. 
The (unique) case with $\dim(\g) = 2$ is also clear, so we will suppose that $\dim(\g) > 2$, and prove that $\h$ is maximal totally isotropic, 
\textit{i.e.} that if $v \in \g$ is a homogeneous element satisfying that $\lambda([v,\h]) = 0$, then $v \in \h$. 
By Remark \ref{rem:superpol} and the fact that $k y$ is an ideal of the super Lie algebra $\k$, we see that $y \in \h$. 
Thus, the assumption that $\lambda([v,\h])=0$ yields that $\lambda([v , y])=0$. 
Using that $y \in \Z(\g/\Z(\g))$ (so $[\g,y] = k z$) and $\lambda(z)=1$, we conclude that the previous vanishing identity is equivalent to $[v,y]=0$, 
\textit{i.e.} $v \in \k$, which in turn implies that $v \in \h$.
\qed

\begin{rem}
\label{rem:lambday}
Note that we can further suppose in the lemma that $\lambda(y)=0$, for $\lambda(y)=0$ if $|y|=1$, and we may change $y$ by $y - \lambda(y) z$ when $y$ is even. 
\end{rem}

Using the lemma we can prove the result:
\begin{teo}
\label{teo:teo1}
Let $\g$ be a nilpotent super Lie algebra, and $\lambda \in \L_{\g}$ be a functional. 
Then there exists a polarization $\h$ of $\g$ at $\lambda$ such that the induced module $\ind(F_{\lambda,\h,i},\g)$ is simple, for $i \in \ZZ/2\ZZ$. 
It may be even assumed that $\h$ is invariant under the action of $\g$. 
\end{teo}
\noindent\textbf{Proof.}
The proof is a variation of that given in \cite{Dix96}, Thm. 6.1.1, but we avoid the use of the so-called standard polarizations. 

We first note that if the super Lie algebra is of dimension at most $2$ and concentrated in one degree, 
\textit{i.e.} $\g=\g_{0}$ or $\g=\g_{1}$, then the theorem is immediate: 
in both cases there is only one polarization $\h = \g$, so 
$\ind(F_{\lambda,\h,i},\g)$ is one-dimensional and the statement also holds. 

We shall now proceed to prove the theorem by induction on the dimension of the underlying vector space of $\g$. 
If $\dim(\g) = 1$ the result follows from the previous considerations. 
If $\dim(\g) = 2$, the only case that does not follow from the previous paragraph is when $\dim(\g_{0})=\dim(\g_{1})=1$. 
Let us suppose that $\g_{0} = k.z$ and $\g_{1}=k.c$. 
By Lemma \ref{lema:abe}, we see that, up to isomorphism, we have two possibilities: $[c,c] = 0$ or $[c,c] = z$. 
Either if we consider the first case for arbitrary $\lambda$ or the second case for $\lambda = 0$, 
there is a unique polarization $\h = \g$, so the theorem holds, for $\ind(F_{\lambda,\h,i},\g)$ is one-dimensional. 
If we regard the second case with $\lambda \neq 0$, we see that $\lambda \in L_{\g}$ and there is a unique polarization $\h = \g_{0}$, which is invariant under the action of $\g$. 
The theorem also holds in this case, because it is a particular case of \cite{Ser99}, Cor. 3.2. 

Let us suppose that $\dim(\g) = d > 2$ and that the proposition holds for dimensions (strictly) less than $d$. 
By Remark \ref{rem:superpol}, we assume there are no ideals of $\g$ such that $\lambda$ vanishes on them. 
In particular, we see that $\Z(\g) \cap \Ker(\lambda)$ should be trivial. 
This implies that $\Z(\g)$ should be one-dimensional, because $\g$ is nilpotent, and included in $\g_{0}$. 
Let $z \in \Z(\g)$ be a nonzero element such that $\lambda(z)=1$. 
By the previous lemma, there exists $x, y \in \g$ homogeneous of the same degree satisfying that $[x,y]=z$, $[y,y]=0$ and $\k=\C(\{y\})$ 
is an ideal of $\g$ such that $\g = \k \oplus k x$. 
Consider $\a = k z \oplus k y$. 
It is clearly a supercommutative ideal of $\g$. 
We see that $\k = \st(\lambda|_{\a},\g)$.  

Let us now suppose that we have a polarization $\h$ of $\k$ at $\lambda|_{\k}$ such that it is invariant under the action of $\g$. 
Such a polarization exists by Remark \ref{rem:polinv}, and by Remark \ref{rem:superpol} it must contain $\a$. 
Lemma \ref{lema:mio} tells us that $\h$ is also a polarization of $\g$ at $\lambda$. 


We have thus a polarization $\h$ of $\g$ at $\lambda$ invariant under the action of $\g$, included in $\k = \st(\lambda|_{\a},\g)$ and including $\a$. 
By the inductive hypothesis $W = \ind(\lambda|_{\h},\k)$ is simple, with structure morphism denoted by $\rho$. 
We remark that the representation $\ind(\lambda|_{\h},\g)$ is obviously isomorphic to $\ind(W,\g)$. 
Using Proposition \ref{prop:5.1.13} we see that the $\a$-representation on $W$ given by $\rho|_{\a}$ is a direct sum of copies of the $\a$-representation given by 
$\lambda|_{\a}$. 
Applying Theorem \ref{teo:5.3.6}, our theorem follows. 
\qed

The next result is a superized version of a lemma appearing in \cite{Dix96}, whose proof applies to this case as well. 
\begin{lema}
\label{lema:subcom}
Let $\a$ be a supercommutative ideal of $\g$ with super centralizer $\C(\a)$, and $\h$ a subalgebra of the super Lie algebra $\g$ subordinate to $\lambda$. 
Set $\bar{\h} = (\h \cap \a^{\lambda}) + \a$, where $\a^{\lambda} = \{ x \in \g : \lambda([x,y])=0, \forall y \in \a \}$. 
Then $\bar{\h}$ is a subalgebra of $\g$ subordinate to $\lambda$ and $\h \cap \C(\a) \cap \ker(\lambda)$ is an ideal of the super Lie algebra $\h + \a$. 
\end{lema}
\noindent\textbf{Proof.}
The proof given in \cite{Dix96}, Lemma 6.1.3, applies word for word.
\qed

We also have the following theorem, whose proof is an adaptation of that in \cite{Dix96}, Thm. 6.1.4.
\begin{teo}
\label{teo:teo2}
Let $\g$ be a nilpotent super Lie algebra and $\lambda \in \Hom(\g,k)$. 
Given two polarizations $\h$ and $\h'$ of $\g$ at $\lambda$, let $\rho_{\h,i}$ and $\rho_{\h',j}$ be the structure morphisms of the $\U(\g)$-modules 
$\ind(F_{\lambda,\h,i},\g)$ and $\ind(F_{\lambda,\h',j},\g)$ determined by the polarizations $\h$ and $\h'$ and by some $i, j \in \ZZ/2 \ZZ$, resp. 
Then, $\ker(\rho_{\h,i}) = \ker(\rho_{\h',j})$. 
\end{teo}
\noindent\textbf{Proof.}
Is is easy to see that $\ker(\rho_{\h,i}) = \ker(\rho_{\h,j})$, for all $i,j \in \ZZ/2\ZZ$, so from now on we will omit the indices $i$ and $j$. 

We first note that if $\g=\g_{0}$ or $\g=\g_{1}$, then the theorem is immediate: 
the first case is just the classical result for Lie algebras (see \cite{Dix96}, Thm. 6.1.4), and in the second case there is only one polarization $\h = \g$, 
so the statement of the theorem holds in both cases. 

We shall now proceed by induction on the dimension of $\g$. 
For $\dim(\g) = 1$ the result is a consequence of the previous considerations. 
If $\dim(\g) = 2$, the only case that does not follow from the previous paragraph is when $\dim(\g_{0})=\dim(\g_{1})=1$. 
Let us suppose that $\g_{0} = k.z$ and $\g_{1}=k.c$. 
By Lemma \ref{lema:abe}, we see that, up to isomorphism, we have two possibilities: $[c,c] = 0$ or $[c,c] = z$. 
Either in the first case for arbitrary $\lambda$ or in the second case for $\lambda = 0$, there is a unique polarization $\h = \g$, for which the theorem holds. 
In the second case with $\lambda \neq 0$, there is also a unique polarization $\h = \g_{0}$, and the statement also follows in this case. 

Let us assume that $\dim(\g) = d > 2$ and that the statement holds for dimensions strictly less than $d$. 
Let $\h$ and $\h'$ be two polarizations of $\g$ at $\lambda$. 

If there exists a nonzero ideal $\k$ such that $\lambda(\k) = 0$, then $\h$ and $\h'$ include $\k$ by Remark \ref{rem:superpol}. 
Passing to the quotient $\g/\k$, we see that $\h/\k$ and $\h'/\k$ are polarizations of $\g/\k$ at the functional $\bar{\lambda}$ induced by $\lambda$. 
Indeed, $\h/\k$ and $\h'/\k$ are obviously subordinate to $\bar{\lambda}$ and maximal. 
The theorem follows in this case by inductive hypothesis. 

We thus suppose that there is no nonzero ideal of $\g$ such that $\lambda(\k) = 0$. 
Since $\g$ is nilpotent, $\Z(\g)$ is a nonzero ideal, so $\dim(\Z(\g))=1$, $\Z(\g) \subseteq \g_{0}$ and $\lambda(\Z(\g)) \neq 0$. 
Set $\Z(\g) = k.z$. 
By Lemma \ref{lema:bm}, there exists $x, y \in \g$ homogeneous of the same degree satisfying that $[x,y]=z$, $[y,y]=0$ 
and $\k=\C(\{y\})$ is an ideal of $\g$ such that $\g = \k \oplus k x$. 
Consider $\a = k z \oplus k y$. 
It is clearly a supercommutative ideal of $\g$. 
Notice that $\k = \st(\lambda|_{\a},\g)$. 

Put $\bar{\h} = (\h \cap \a^{\lambda}) + \a$ and $\bar{\h}' = (\h' \cap \a^{\lambda}) + \a$. 
Lemma \ref{lema:subcom} tells us that $\bar{\h}$ and $\bar{\h}'$ are subordinate to $\lambda$. 
Both of them satisfy that $\bar{\h}, \bar{\h}' \subseteq \k$ and also $\sdim(\h) = \sdim(\bar{\h})$ and $\sdim(\h') = \sdim(\bar{\h}')$. 
If we restrict to $\k$, $\bar{\h}$ and $\bar{\h}'$ are subalgebras subordinated to $\lambda|_{\k}$. 
Furthermore, they are polarizations of $\k$ at $\lambda|_{\k}$, because they have the same super dimension as $\h$ and $\h'$, resp. 
By inductive hypothesis, they satisfy that $\ker(\ind(\lambda|_{\bar{\h}},\k)) = \ker(\ind(\lambda|_{\bar{\h}'},\k))$, so by Proposition \ref{prop:5.1.7}, 
we have that $\ker(\ind(\lambda|_{\bar{\h}},\g)) = \ker(\ind(\lambda|_{\bar{\h}'},\g))$. 

We must then show that $\ker(\rho_{\bar{\h}}) = \ker(\rho_{\h})$ in order to conclude the proof. 
Since, if $\h \subseteq \k$, then $\h = \bar{\h}$ (because any polarization $\h$ included in $\k$ should satisfy that $\a \subseteq \h$, for $\lambda([\a,\k])=0$), 
we shall assume that $\h \nsubseteq \k$. 
In this case, we may further suppose that $x \in \h$ (by the proof of \cite{BM90}, Lemma 1.10). 
Set $\n = \h + \a$. 
Since $\lambda([x,y])= \lambda(z) \neq 0$, we see that $y \notin \h$. 
Also note that $z \in \h$ (by Remark \ref{rem:superpol}). 
We see that $\h, \bar{\h} \subseteq \n$ are polarizations of $\n$ at $\lambda|_{\n}$, because they are subordinated to $\lambda|_{\n}$ and of the appropriate super dimension. 
By Proposition \ref{prop:5.1.7} we see that it suffices to prove that $\ker(\ind(\lambda|_{\h},\n)) = \ker(\ind(\lambda|_{\bar{\h}},\n))$. 
Since, by Lemma \ref{lema:subcom}, $\h \cap \k \cap \ker(\lambda)$ is an ideal in $\n$, by inductive hypothesis we will suppose that the former is trivial. 
Then, $\dim(\h \cap \k) \leq 1$, so $\h \cap \k = \Z(\g)$, and analogously for $\bar{\h}$. 
This implies that 
\begin{align*}
 \n &= k z \oplus k y \oplus k x,
 \\
 \h &= k z \oplus k x,
 \\
 \bar{\h} &= k z \oplus k y.
\end{align*}
If $|x| = |y| = 0$, the statement follows from \cite{Dix96}, Lemma 6.1.2, (iii). 
If $|x| = |y| = 1$, the statement follows from \cite{Ser99}, Lemma 1.2, 2). 
The theorem is thus proved. 
\qed

From Theorem \ref{teo:teo1} we see that given $\lambda \in \L_{\g}$, there exist a primitive ideal $I(\lambda)$ of $\U(\g)$ given as the kernel of the structure morphism of the representation $\ind(\lambda|_{\h},\g)$, for some polarization $\h$ of $\g$ at $\lambda$. 
In fact, by Theorem \ref{teo:teo2}, the ideal does not depend on the polarization. 

The following proposition follows from the work of Letzter and the previous theorems, and in fact provides a link between our point of view and his. 
\begin{prop}
\label{prop:let}
Let $\lambda \in \L_{\g}$ be a functional of a nilpotent super Lie algebra $\g$. 
Then, there exists a unique maximal ideal of $\U(\g_{0})$ containing $I(\lambda) \cap \U(\g_{0})$, and it is in fact $I(\lambda|_{\g_{0}})$. 
\end{prop}
\noindent\textbf{Proof.}
Let $\lambda_{0} \in \g_{0}^{*}$ be the restriction of $\lambda$ to $\g_{0}$, and $I(\lambda_{0})$ be primitive ideal of $\U(\g_{0})$ determined by it. 
Set $\h$ a polarization of $\g$ at $\lambda$. 
Then $\h_{0}$ is a polarization of $\g_{0}$ at $\lambda_{0}$. 
If we also denote by $\lambda$ and $\lambda_{0}$ the structure morphisms of the one-dimensional representations over $\U(\h)$ and $\U(\h_{0})$ that they determine, respectively, we have the commutative diagram 
\[
\xymatrix
{
\U(\h)
\ar[drr]^-{\lambda}
&
&
\\
&
&
k
\\
\U(\h_{0})
\ar[urr]^-{\lambda_{0}}
\ar@{^{(}->}[uu]
} 
\]
If $J$ is the kernel of $\lambda$ and $J_{0}$ the kernel of $\lambda_{0}$, the commutativity of the diagram says that $J_{0} \subseteq J$. 
Moreover, the PBW Theorem tells us that $J \cap \U(\h_{0}) = J_{0}$ and also that $\U(\g) J \cap \U(\g_{0}) = \U(\g_{0}) J_{0}$. 
By construction, $I(\lambda_{0})$ is the largest ideal of $\U(\g_{0})$ inside of $\U(\g_{0}) J_{0}$ and the ideal 
$I(\lambda)$ is the largest ideal of the super algebra $\U(\g)$ inside of $\U(\g) J$. 
Hence, $I(\lambda) \cap \U(\g_{0}) \subseteq \U(\g) J \cap \U(\g_{0}) = \U(\g_{0}) J_{0}$, which yields that $I(\lambda) \cap \U(\g_{0}) \subseteq I(\lambda_{0})$. 
Since, by \cite{Let92}, Cor. III, we have that $I(\lambda) \cap \U(\g_{0})$ has a unique minimal prime ideal, which is a primitive ideal of $\U(\g_{0})$, 
the previous inclusion implies that it must be $I(\lambda_{0})$. 
\qed

The following result tells us that every primitive ideal is of the form $I(\lambda)$, for some functional $\lambda \in \L_{\g}$ 
(\textit{cf.}~\cite{Dix96}, Thm. 6.1.7). 
\begin{teo}
\label{teo:teo3}
Let $I$ be a primitive ideal of the enveloping algebra $\U(\g)$ of a nilpotent super Lie algebra $\g$. 
Then, there exists $\lambda \in \L_{\g}$ such that $I=I(\lambda)$. 
\end{teo}
\noindent\textbf{Proof.} 
We may derive this theorem as a consequence of the work of Letzter. 
Since $I$ is a primitive ideal of $\U(\g)$, then \cite{Let92}, Cor. III, yields that $I \cap \U(\g_{0})$ has a unique minimal prime ideal, 
which is a primitive ideal of $\U(\g_{0})$, and this assignment is in fact a bijection. 
Let $\lambda_{0} \in \g_{0}^{*}$ be a linear functional such that $I(\lambda_{0})$ is the previous primitive ideal of $\U(\g_{0})$, and 
$\lambda \in \L_{\g}$ the obvious extension of $\lambda_{0}$ to $\g$, so $\lambda_{0} = \lambda|_{\g_{0}}$. 
The previous proposition tells us that $I(\lambda) \cap \U(\g_{0}) \subseteq I(\lambda_{0})$ and, by construction, we have that 
$I \cap \U(\g_{0}) \subseteq I(\lambda_{0})$. 
Since the map given in \cite{Let92}, Cor. III, is bijective, we conclude that $I = I(\lambda)$. 
The theorem is thus proved. 
\qed

\begin{lema}
\label{lema:bm2}
Let $\g$ be a nilpotent super Lie algebra with super center $\z = k z \neq \g$, with $z$ even, $x, y \in \g$ homogeneous of the same parity and $\k$ an ideal 
of codimension one in $\g$ satisfying the properties (i), (ii) and (iii) and (2) stated in Lemma \ref{lema:bm}. 
Define $\bar{\k} = \k/k y$, $\delta$ the locally nilpotent derivation of $\U(\k)$ induced by $\ad (x)$, and $\bar{u}$ the corresponding image element of $u \in \U(\k)$ 
under the canonical projection from $\U(\k)$ to $\U(\bar{\k})$. 
Then,
\begin{itemize}
\item[(i)] if $x$ is even, there exists a unique morphism of super algebras $\phi_{0} : \U(\g) \rightarrow \U(\bar{\k}) \otimes A_{1}(k)$ given by 
\begin{align*}
   \phi_{0} (x) &= 1 \otimes p,
   \\
   \phi_{0} (u) &= \sum_{n \in \NN_{0}} \frac{1}{n!} \overline{\delta^{n}(u)} \otimes q^{n},
\end{align*}
where $u \in \U(\k)$, and $A_{1}(k)$ is the super algebra described in Example \ref{ex:weylcliff}, (i). 
Moreover, it induces an isomorphism $\psi_{0}$ from $\U(\g)_{z}$ to $\U(\bar{\k})_{z} \otimes A_{1}(k)$.

\item[(ii)] if $x$ is odd, there exists a unique morphism of super algebras $\phi_{1} : \U(\g) \rightarrow \U(\bar{\k}) \otimes M_{2}(k)$ given by 
\begin{align*}
   \phi_{1} (x) &= \begin{pmatrix} 0 & \frac{\overline{[x,x]}}{2} \\ 1 & 0 \end{pmatrix},
   \\
   \phi_{1} (u) &= \begin{pmatrix} \bar{u} & \overline{\delta(u)} \\ 0 & \overline{\Sigma(u)} \end{pmatrix},
\end{align*}
where $u \in \U(\k)$, and $M_{2}(k)$ is the super algebra described in Example \ref{ex:weylcliff}, (ii). 
Furthermore, it induces an isomorphism $\psi_{1}$ from $\U(\g)_{z}$ to $\U(\bar{\k})_{z} \otimes M_{2}(k)$. 
\end{itemize} 
Moreover, given $I \neq \U(\g)$ an ideal of $\U(\g)$ such that $z-1 \in I$, then there exists one and only one ideal 
$J$ of $\U(\bar{\k})$ satisfying that $\bar{z}-1 \in J$ and $\phi_{0}(I_{z}) = J_{\bar{z}} \otimes A_{1}(k)$, if $x$ is even, or 
$\phi_{1}(I_{z}) = J_{\bar{z}} \otimes M_{2}(k)$, if $x$ is odd.  
Finally, there is a chain of isomorphisms of super algebras 
\[     \U(\g)/I \rightarrow \U(\g)_{z}/I_{z} \rightarrow (\U(\bar{\k})_{\bar{z}}/J_{\bar{z}}) \otimes A_{1}(k) \rightarrow (\U(\bar{\k})/J) \otimes A_{1}(k),     \]
if $x$ is even, and 
\[     \U(\g)/I \rightarrow \U(\g)_{z}/I_{z} \rightarrow (\U(\bar{\k})_{\bar{z}}/J_{\bar{z}}) \otimes M_{2}(k) \rightarrow (\U(\bar{\k})/J) \otimes M_{2}(k),     \]
if $x$ is odd.
\end{lema}
\noindent\textbf{Proof.}
The proof of the first statement of both items is implicit in \cite{BM90}, and follows the lines of \cite{Dix96}, Lemmas 4.6.6 and 4.7.8, (i), but we give it for clarity. 

\cite{BM90}, Lemma 1.7, tells us that there exist isomorphisms of super algebras $\U(\g) \simeq \U(\h)[t,\id,\delta]$, for $x$ even, 
given by $u \mapsto u$, if $u \in \U(\h)$, and $x \mapsto t$, 
and $\U(\g) \simeq \U(\h)[t,\Sigma,\delta]/(t^{2}-[x,x]/2)$, for $x$ odd, given by $u \mapsto \bar{u}$, if $u \in \U(\h)$, and $x \mapsto \bar{t}$, 
where the classes here are with respect to the quotient by the ideal $(t^{2} - [x,x]/2)$. 
They obviously induce isomorphisms $\U(\g)_{z} \simeq \U(\h)_{z}[t,\id,\delta]$, for $x$ even, and 
$\U(\g)_{z} \simeq \U(\h)_{z}[t,\Sigma,\delta]/(t^{2}-[x,x]/2)$, for $x$ odd, respectively. 
Furthermore, \cite{BM90}, Lemmas 1.4 and 1.5, gives explicit isomorphisms $\U(\h)_{z}[t,\id,\delta] \simeq \U(\bar{\h})_{\bar{z}} \otimes A_{1}(k)$, 
for $x$ even, and $\U(\h)_{z}[t,\Sigma,\delta]/(t^{2}-[x,x]/2) \simeq \U(\bar{\h})_{\bar{z}} \otimes M_{2}(k)$, for $x$ odd, 
which are just the maps $\psi_{0}$ and $\psi_{1}$ in the statement, respectively. 
Making the composition of the previous morphisms with the canonical map $\U(\g) \rightarrow \U(\g)_{z}$, we obtain maps 
$\U(\g) \rightarrow \U(\bar{\h})_{\bar{z}} \otimes A_{1}(k)$, for $x$ even, and $\U(\g) \rightarrow \U(\bar{\h})_{\bar{z}} \otimes M_{2}(k)$, for $x$ odd, 
respectively. 
It is trivial to check that the images of these morphisms are in fact contained in the image of the canonical maps 
$\U(\bar{\h}) \otimes A_{1}(k) \rightarrow \U(\bar{\h})_{\bar{z}} \otimes A_{1}(k)$ and 
$\U(\bar{\h}) \otimes M_{2}(k) \rightarrow \U(\bar{\h})_{\bar{z}} \otimes M_{2}(k)$, respectively, giving us the morphisms $\phi_{0}$ and $\phi_{1}$, 
respectively. 

Finally, the proof of the two chain of isomorphisms on item (ii) is the same as the one given in \cite{Dix96}, Lemma 4.7.8, (ii), 
replacing the use of \cite{Dix96}, Proposition 3.6.15 and Lemma 4.5.1 by Proposition \ref{prop:3.6.15} and Lemma \ref{lema:4.5.1}, respectively. 
\qed

\begin{lema}
\label{lema:6.2.1}
Let $\g$ be a super Lie algebra with super center $\z = k z \neq \g$, where $z$ is even, $x, y \in \g$ be homogeneous elements of the same parity 
and $\k$ an ideal of codimension one in $\g$ satisfying the properties (i), (ii), (iii) and (2) stated in Lemma \ref{lema:bm}. 
Set $\bar{\k} = \k/k.y$. 
Given $\lambda \in \L_{\g}$ satisfying that $\lambda(z) = 1$ and $\lambda(y)=0$, 
define $\lambda' = \lambda|_{\k} \in \L_{\k}$ and $\bar{\lambda}'$ the functional induced by $\lambda'$ on $\bar{\k}$. 
Then, $\psi_{0}(I(\lambda)_{z}) = I(\bar{\lambda}')_{\bar{z}} \otimes A_{1}(k)$ if $x$ is even, and $\psi_{1}(I(\lambda)_{z}) = I(\bar{\lambda}')_{\bar{z}} \otimes M_{2}(k)$ 
if $x$ is odd. 
\end{lema}
\noindent\textbf{Proof.}
The proof for $x$ even is the same as the one appearing in \cite{Dix96}, Lemma 6.2.1, with the additional assumption that all elements must be homogeneous. 
The proof for $x$ odd is exactly the same as in the even case, replacing the use of \cite{Dix96}, Lemma 4.7.8 and Prop. 5.1.7, by Lemma \ref{lema:bm2} and Proposition \ref{prop:5.1.7}, 
respectively, and the appearances of the Weyl algebra $A_{1}(k)$ by $M_{2}(k)$ and of $1 \otimes q$ by $1 \otimes e_{12}$, and using that $M_{2}(k)$ is also a simple super algebra. 
\qed

The final result of this section is the following. 
\begin{prop}
\label{prop:6.2.3}
Let $\g$ be a nilpotent super Lie algebra and let $\mathcal{A}d_{0}$ be the adjoint group of the Lie algebra $\g_{0}$, acting on $\g^{*}_{0}$. 
Given $\lambda, \lambda' \in \L_{\g}$, then $I(\lambda) = I(\lambda')$ if and only if $\lambda'$ and $\lambda$ lie in the same orbit of $\g^{*}_{0}$ 
under the coadjoint action of $\mathcal{A}d_{0}$. 
\end{prop}
\noindent\textbf{Proof.}
Even though a proof following the lines of \cite{Dix96}, Prop. 6.2.3, is possible, we give a shorter one. 
 
For $a \in \mathcal{A}d_{0}$, we denote by $a_{\U}$ the automorphism of $\U(\g)$ induced by $a$. 
Suppose that $a(\lambda) = \lambda'$. 
Then, by transport of structures $a_{\U}(I(\lambda)) = I(\lambda')$, and using Proposition \ref{prop:2.4.17} we conclude that $I(\lambda) = I(\lambda')$. 

Conversely, let us assume that $I(\lambda) = I(\lambda')$. 
Then $I(\lambda) \cap \U(\g_{0}) = I(\lambda') \cap \U(\g_{0})$, so, by Proposition \ref{prop:let}, we have that $I(\lambda|_{\g_{0}}) = I(\lambda'|_{\g_{0}})$. 
Now, \cite{Dix96}, Prop. 6.2.3, tells us that there exist $a \in \mathcal{A}d_{0}$ such that $a(\lambda|_{\g_{0}}) = \lambda'|_{\g_{0}}$, 
which \textit{a fortiori} yields that $a(\lambda) = \lambda'$. 
The proposition is thus proved. 
\qed

\subsection{Some consequences}
\label{subsec:sevcons}

We want to derive some consequences from the main theorems proved before. 

\subsubsection{Simple quotients of the enveloping algebra of a nilpotent super Lie algebra}

From Lemma \ref{lema:6.2.1} we obtain the following proposition, which is analogous to \cite{Dix96}, Prop. 6.2.2 (\textit{cf.}~\cite{BM90}, Thm. A).
\begin{prop}
\label{prop:mio}
Let $\g$ be a nilpotent super Lie algebra and $\lambda \in \L_{\g}$. 
The primitive ideal $I(\lambda)$ satisfies that $\U(\g)/I(\lambda) \simeq \Cliff_{q}(k) \otimes A_{p}(k)$, where 
$(p,q) = (\dim(\g_{0}/\g_{0}^{\lambda})/2,\dim(\g_{1}/\g_{1}^{\lambda}))$ and $\g^{\lambda} = (\g_{0}^{\lambda},\g_{1}^{\lambda})$ is the kernel of 
the superantisymmetric bilinear form $\cl{\hskip 0.6mm, }_{\lambda}$ determined by $\lambda$ on $\g$. 
\end{prop}
\noindent\textbf{Proof.}
We first remark that if $\g=\g_{0}$ or $\g=\g_{1}$, then the proposition is immediate. 
Indeed, the first case is just the classical result for Lie algebras (see \cite{Dix96}, Prop. 6.2.2). 
In the second case, $\g^{\lambda} = \g$ and $\U(\g)$ is a supersymmetric super algebra with a unique maximal ideal $I$ 
whose quotient is $k$. 

We shall now proceed by induction on the dimension of $\g$. 
If $\dim(\g)=1$, the result follows from the previous considerations. 
In case $\dim(\g) = 2$, the only case that does not follow from the previous paragraph is when $\dim(\g_{0})=\dim(\g_{1})=1$. 
Let us suppose that $\g_{0} = k.z$ and $\g_{1}=k.c$. 
By Lemma \ref{lema:abe}, we have two possibilities up to isomorphism: $[c,c] = 0$ or $[c,c] = z$. 
It is not difficult to prove that, either if we consider the first case for arbitrary $\lambda$ or the second case for $\lambda = 0$, 
$\g^{\lambda} = \g$, so the theorem holds, for $\ind(F_{\lambda,\h,i},\g)$ is one-dimensional. 
If we regard the second case with $\lambda \neq 0$, we see that $\g^{\lambda} = \g_{0}$. 
It can be easily checked that the annihilator is the ideal generated by $z-1$, whose quotient is $\Cliff_{1}(k)$, so the statement also holds in this case 
(\textit{cf.}~\cite{BM90}, 0.2, (b)). 

Let us suppose that $\dim(\g) > 2$ and let $\z$ denote the super center of $\g$. 
We denote $I(\lambda)$ simply by $I$ and consider $\h$ a polarization of $\g$ at $\lambda$ such that $\ind(\lambda|_{\h},\g)$ is simple. 

If $I \cap \z \neq 0$, then $\lambda(I \cap \z)=0$ and $I \cap \z \subseteq \h$, so we may consider $\bar{\lambda} \in \L_{\g/(I \cap \z)}$ induced by $\lambda$. 
It is easy to see that $\bar{\h} = \h/(I \cap \z)$ is a polarization at $\bar{\lambda}$ 
and that $\ind(\bar{\lambda}|_{\bar{\h}},\g/(I \cap \z))$ is a simple $\U(\g/(I \cap \z))$-module. 
Moreover, the image of $I$ under the projection map $\U(\g) \rightarrow \U(\g/(I \cap \z))$ coincides with the kernel $\bar{I}$ of the structure morphism of 
$\ind(\bar{\lambda}|_{\bar{\h}},\g/(I \cap \z))$, thus $\U(\g)/I \simeq \U(\g/(I \cap \z))/\bar{I}$. 
It is also clear that $\g^{\lambda} \supseteq I \cap \z$, so $(\g/(I \cap \z))^{\bar{\lambda}} = \g^{\lambda}/(I \cap \z)$ and 
$(\g/(I \cap \z))/(\g/(I \cap \z))^{\lambda} = \g/\g^{\lambda}$. 
Then, the statement follows from the inductive hypothesis. 

Let us now assume that $I \cap \z = 0$, which tells us that $\dim(\z)=1$. 
Suppose that $\z = k z \neq \g$, where $z$ is even, and consider $x, y \in \g$ be homogeneous of the same parity and $\k$ an ideal 
of codimension one in $\g$ satisfying the properties (i), (ii), (iii) and (2) stated in Lemma \ref{lema:bm}, such that $\lambda(z) = 1$ and $\lambda(y)=0$, so $z - 1 \in I$. 
Set $\bar{\k} = \k/k.y$ and define $\lambda' = \lambda|_{\k} \in \L_{\k}$ and $\bar{\lambda}'$ the functional induced by $\lambda'$ on $\bar{\k}$. 
It is direct to check that $\g^{\lambda} \subseteq \k$ and moreover $\sdim(\bar{\k}^{\bar{\lambda}'}) = \sdim(\g^{\lambda})$. 
By Lemma \ref{lema:6.2.1}, $\psi_{0}(I(\lambda)_{z}) = I(\bar{\lambda}')_{\bar{z}} \otimes A_{1}(k)$ if $x$ is even, and $\psi_{1}(I(\lambda)_{z}) = I(\bar{\lambda}')_{\bar{z}} \otimes M_{2}(k)$ if $x$ is odd. 
The corollary thus follows from the inductive assumption and Lemma \ref{lema:bm2}. 
\qed

\subsubsection{Maximal ideals of the underlying algebra of the enveloping algebra of a nilpotent super Lie algebra}

This paragraph is devoted to obtain a ``parametrization'' of the maximal ideals of the underlying algebra of $\U(\g)$, similar 
to the one given for maximal ideals. 

By Lemma \ref{lema:primax} we know that, for every maximal ideal $I$ of a super algebra $A$, the set of minimal prime ideals $J$ of the underlying algebra $\O(A)$ of $A$ 
such that $J \supseteq I$ form an orbit under $\Sigma$, which \textit{a fortiori} has at most two elements, and they are in fact maximal ideals of the underlying algebra. 
Moreover, all the maximal ideals of $\O(A)$ can be realized in this way, and given any $J$ as before, we have that $I = J \cap \Sigma(J)$. 

The following proposition gives us a description of the maximal ideals of the underlying algebra of $\U(\g)$, more or less parallel to the one given for 
the construction of the ideals $I(\lambda)$. 
\begin{prop}
Let $\g$ be a nilpotent super Lie algebra, $\lambda \in \L_{\g}$ a functional and 
$I(\lambda)$ be the corresponding maximal ideal of the super algebra $\U(\g)$. 
If $\dim(\g_{1}/\g_{1}^{\lambda})$ is even, then $I(\lambda)$ is a maximal ideal of the underlying algebras of $\U(\g)$. 
If $\dim(\g_{1}/\g_{1}^{\lambda})$ is odd, then $I(\lambda)$ is not a maximal ideal of $\O(\U(\g))$, 
and there exist two maximal ideals $I^{+}(\lambda)$ and $I^{-}(\lambda)$ of the underlying algebra of $\U(\g)$, such that 
$\Sigma(I^{-}(\lambda))=I^{+}(\lambda)$ and $I^{+}(\lambda) \cap I^{-}(\lambda) = I(\lambda)$. 
Furthermore, any maximal ideal containing $I(\lambda)$ is one of these, and all the maximal ideals of $\O(\U(\g))$ 
can be obtained in this way. 

More explicitly, the ideals $I^{+}(\lambda)$ and $I^{-}(\lambda)$ can be realized as follows. 
Given a polarization $\h$ of $\g$ at $\lambda$ invariant under the action of $\g$ such that $\ind(\lambda|_{\h},\g)$ is a simple module 
over the super algebra $\U(\g)$, there exists an element $c \in \g_{1}$ such that $\lambda([c,c]) = 2$ and $\lambda([c,\h]) = 0$, 
which also satisfies that $[c,c] \in \h$. 
Define the subalgebra $\hat{\h} = \h \oplus k.c$ of $\g$. 
Then there exists two extensions $\lambda_{+}$ and $\lambda_{-}$ of the structure morphism $\lambda : \U(\h) \rightarrow k$ 
of the one-dimensional $\h$-representation given by $F_{\lambda,\h}$ to morphisms of algebras $\lambda_{\pm} : \U(\hat{\h}) \rightarrow k$. 
We denote these modules by $F_{\lambda,\h,\pm}$. 
If $J$ is the kernel of $\lambda$ and $J^{\pm}$ the kernel of $\lambda_{\pm}$, then $\U(\g) J^{+} \cap \U(\g) J^{-} = \U(\g) J$. 
Finally, $I^{\pm}(\lambda)$ is the largest ideal of $\O(\U(\g))$ inside of $\U(\g) J^{\pm}$.
\end{prop}
\noindent\textbf{Proof.}
By the comments previous to the proposition, we have two possibilities: either $I(\lambda)$ is also a maximal ideal of the underlying algebra of $\U(\g)$, 
or the orbit of maximal ideals of $\O(\U(\g))$ containing $I(\lambda)$ has two elements $J$ and $\Sigma(J)$. 
It is easy to see that the first possibility occurs exactly if $\dim(\g_{1}/\g_{1}^{\lambda})$ is even and the second one when 
$\dim(\g_{1}/\g_{1}^{\lambda})$ is odd. 
This can be proved as follows (\textit{cf.}~\cite{BM90}). 
By the Proposition \ref{prop:mio}, there is an isomorphism of super algebras (and hence also an isomorphism of the underlying algebras) 
$\U(\g)/I(\lambda) \simeq \Cliff_{q}(k) \otimes A_{p}(k)$, where $(p,q) = (\dim(\g_{0}/\g_{0}^{\lambda})/2,\dim(\g_{1}/\g_{1}^{\lambda}))$. 
If $q = \dim(\g_{1}/\g_{1}^{\lambda})$ is even, then $\Cliff_{q}(k) \simeq M_{2^{q/2}}(k)$, so $\U(\g)/I(\lambda)$ is a simple algebra, 
which tells us that $I(\lambda)$ is a maximal ideal of the underlying algebra of $\U(\g)$. 
On the other hand, if $q = \dim(\g_{1}/\g_{1}^{\lambda})$ is odd, then $\Cliff_{q}(k) \simeq M_{2^{(q-1)/2}}(k) \otimes k[\epsilon]/(\epsilon^{2} - 1)$. 
Using the obvious isomorphism of algebras $k[\epsilon]/(\epsilon^{2}-1) \simeq k \times k$, given by $a + b \epsilon \mapsto (a + b, a - b)$, 
for $a, b \in k$, we obtain the isomorphisms of algebras $\Cliff_{q}(k) \simeq M_{2^{(q-1)/2}}(k) \times M_{2^{(q-1)/2}}(k)$ 
and, in consequence, $\U(\g)/I(\lambda) \simeq M_{2^{(q-1)/2}}(A_{p}(k)) \times M_{2^{(q-1)/2}}(A_{p}(k))$. 
We conclude that in this case $I(\lambda)$ is not a maximal ideal of the underlying algebra of $\U(\g)$, and any of the two maximal ideal $J$ 
of $\O(\U(\g))$ containing $I(\lambda)$ satisfy that $\U(\g)/J \simeq M_{2^{(q-1)/2}}(A_{p}(k))$. 

If $q = \dim(\g_{1}/\g_{1}^{\lambda})$ is even, $I(\lambda)$ is already a maximal ideal of $\O(\U(\g))$, so there is nothing to do. 
Let us assume that it is odd and we write $q = 2 r + 1$. 
We denote $s$ the dimension of $\dim(\g_{1}^{\lambda})$. 
This implies that there exists a basis $\{ y_{1}, \dots, y_{q}, z_{1}, \dots, z_{s} \}$ of $\g_{1}$, for $\{ z_{1}, \dots, z_{s} \}$ a basis of $\g_{1}^{\lambda}$, 
such that the matrix of the symmetric bilinear form $B_{\lambda}$ on that basis is of the form 
\[     \begin{pmatrix} 0_{r \times r} & \mathrm{Id}_{r} & 0_{r \times 1} & 0_{r \times s} 
                       \\
                       \mathrm{Id}_{r} & 0_{r \times r} & 0_{r \times 1} & 0_{r \times s} 
                       \\
                       0_{1 \times r} & 0_{1 \times r} & 2 & 0_{1 \times s} 
                       \\
                       0_{s \times r} & 0_{s \times r} & 0_{s \times 1} & 0_{s \times s}
       \end{pmatrix},     \]
where $0_{m \times n}$ denotes a $m \times n$-matrix with zero entries and $\mathrm{Id}_{r}$ is the identity matrix of $M_{r}(k)$. 
We even choose the previous basis in such a way that $\{ y_{1}, \dots, y_{r}, z_{1}, \dots, z_{s} \}$ is the basis of $\h_{1}$, where $\h$ is a polarization of $\g$ at $\lambda$ 
invariant under the action of $\g$ (\textit{cf.}~Remark \ref{rem:polinv}). 
Fix also $\{ x_{1}, \dots, x_{t} \}$ a basis of $\g_{0}$. 
Note that $\lambda([\h,y_{q}]) = 0$, by construction. 
This implies that $[y_{q},y_{q}] \in \h$. 
Indeed, the fact that $\h$ is invariant under the action of $\g$ tells us that 
$\lambda([h,[y_{q},y_{q}]]) = \lambda([[h,y_{q}],y_{q}]) + (-1)^{|h|} \lambda([y_{q},[h,y_{q}]]) = 0$, for all $h \in \h$. 
Since $\h$ is a polarization, it must be that $[y_{q},y_{q}] \in \h$, as claimed. 

Now, write $c = y_{q}$, set $\hat{\h}$ to be the subalgebra $\h \oplus k.c$ of $\g$, and 
consider the one-dimensional modules $M_{\pm}$ over the underlying algebra of $\U(\hat{\h})$ given by the vector space $F_{\lambda,\h,\pm} = k.v_{\lambda,\h,\pm}$ 
provided with the action $h.v_{\lambda,\h,\pm} = \lambda(h) v_{\lambda,\h,\pm}$, for $h \in \h$, and 
$c.v_{\lambda,\h,\pm} = \pm v_{\lambda,\h,\pm}$. 
A trivial computation shows that it is well-defined and it is in fact an extension of the one-dimensional representation of $\U(\h)$ given by $\lambda$. 
We shall denote the structure morphisms of $F_{\lambda,\h,\pm}$ by $\lambda_{\pm}$. 
We remark that these modules depend on the choice of $c$, but we do not include it in the notation for simplicity. 
Note that $\lambda_{\pm} = \lambda \pm c^{*}$, where $c^{*} = y_{q}^{*}$ is the corresponding functional of the dual basis of 
$\{ x_{1}, \dots, x_{t}, y_{1}, \dots, y_{q}, z_{1}, \dots, z_{s} \}$. 

We have the following commutative diagram of algebras
\[
\xymatrix
{
\U(\hat{\h}) 
\ar[dr]^-{\Sigma}
\ar[drrr]^-{\lambda_{-}}
&
&
&
\\
&
\U(\hat{\h}) 
\ar[rr]^-{\lambda_{+}}
&
&
k
\\
\U(\h)
\ar[dr]^-{\Sigma}
\ar[urrr]^-{\lambda}
\ar@{^{(}->}[uu]
&
&
&
\\
&
\U(\h)
\ar[uurr]^-{\lambda}
\ar@{^{(}->}[uu]
}
\]
Let us denote by $J$ the kernel of $\lambda$ and by $J^{\pm}$ the kernel of $\lambda_{\pm}$. 
We easily see that $\Sigma(J^{-})=J^{+}$ and that $\U(\hat{\h}) J = J^{+} \cap J^{-}$. 
Indeed, from the commutativity of the diagram, we see that $J \subseteq J^{+}$ and $J \subseteq J^{-}$, 
so $\U(\hat{\h}) J \subseteq J^{+} \cap J^{-}$. 
Let us now suppose that $z \in \U(\hat{\h})$ and write $z = u + c u'$,
where $u, u' \in \U(\h)$ are uniquely determined, by the PBW Theorem. 
Then, $z$ belongs to $J^{+} \cap J^{-}$ if and only if 
\begin{align*}
   \lambda_{+}(z) = \lambda(u) + \lambda(u') = 0,
   \\
   \lambda_{-}(z) = \lambda(u) - \lambda(u') = 0.
\end{align*}
Hence, $\lambda(u) = \lambda(u') = 0$, so $z \in \U(\hat{\h}) J$, as we wanted to prove. 
This in turn yields the isomorphism $\U(\hat{\h}) \otimes_{\U(\h)} F_{\lambda,\h} \simeq F_{\lambda,\h,+} \oplus F_{\lambda,\h,-}$ of modules over the 
underlying algebra of $\U(\hat{\h})$. 
Since $\U(\g)$ is a free module over the algebra $\O(\U(\hat{\h}))$, the functor $\U(\g) \otimes_{\U(\hat{\h})} (\place)$ is exact and preserves pull-backs, 
which tells us that $\U(\g) J^{+} \cap \U(\g) J^{-} = \U(\g) J$ and that there exist an isomorphism 
$\ind(\lambda|_{\h},\g) \simeq (\U(\g) \otimes_{\U(\hat{\h})} F_{\lambda,\h,+}) \oplus (\U(\g) \otimes_{\U(\hat{\h})} F_{\lambda,\h,-})$ 
of modules over $\O(\U(\g))$. 
We remark that, by construction, $\U(\hat{\h}) J \subsetneq J^{\pm}$, which further implies that $\U(\g) J \subsetneq \U(\g) J^{\pm}$, by the PBW Theorem. 

We shall now prove the last claim. 
Let $I^{\pm}(\lambda)$ be the largest ideal of $\O(\U(\g))$ inside of $\U(\g) J^{\pm}$. 
It is easy to see, using an argument similar to the proof of \cite{Dix96}, Prop. 5.1.7, that $I^{\pm}(\lambda)$ is the annihilator 
of the module $\U(\g) \otimes_{\U(\hat{\h})} F_{\lambda,\h,\pm}$ over the algebra $\O(\U(\g))$. 
Hence, $\Sigma(I^{+}(\lambda)) \subseteq I^{-}(\lambda)$ and $\Sigma(I^{-}(\lambda)) \subseteq I^{+}$, which imply the equality 
$\Sigma(I^{+}(\lambda)) = I^{-}(\lambda)$. 
By Proposition \ref{prop:5.1.7}, we have that $I(\lambda)$ is the largest ideal of the super algebra $\U(\g)$ inside $\U(\g) J$. 
Since $I^{+}(\lambda) \cap I^{-}(\lambda)$ is an ideal of the super algebra $\U(\g)$ included in $\U(\g) J^{+} \cap \U(\g) J^{-} = \U(\g) J$, then 
$I^{+}(\lambda) \cap I^{-}(\lambda) \subseteq I(\lambda)$. 
Conversely, $I$ is an ideal of $\O(\U(\g))$ included in both $\U(\g) J^{+}$ and $\U(\g) J^{-}$, which yields that 
$I(\lambda) \subseteq I^{+}(\lambda)$ and $I(\lambda) \subseteq I^{-}(\lambda)$, so $I(\lambda) \subseteq I^{+}(\lambda) \cap I^{-}(\lambda)$, 
and thus the equality $I^{+}(\lambda) \cap I^{-}(\lambda) = I(\lambda)$ holds. 

We still have to prove that the ideals $I^{+}(\lambda)$ and $I^{-}(\lambda)$ are maximal. 
Since they are the annihilators of $M_{+} = \U(\g) \otimes_{\U(\hat{\h})} F_{\lambda,\h,+}$ and $M_{-} = \U(\g) \otimes_{\U(\hat{\h})} F_{\lambda,\h,-}$ 
it suffices to show that they are simple modules over the underlying algebra of $\U(\g)$. 

In order to do so, we shall first need the following easy fact. 
As usual, given a simple module $M$ with structure morphism $\rho$ over the underlying algebra of a super algebra $A$, 
we can consider the module $M^{\Sigma}$ over $\O(A)$ with the same underlying vector space but with the structure morphism defined by $\rho \circ \Sigma$. 
If $M$ is a simple module over the underlying algebra $\O(A)$, then $M \oplus M^{\Sigma}$ has the structure of a module over the super algebra $A$. 
This can be proved as follows. 
Let $\mathfrak{m} \subseteq \O(A)$ be a maximal left ideal of the underlying algebra of $A$ such that $A/\mathfrak{m} \simeq M$. 
It is clear that $M^{\Sigma} \simeq A/\Sigma(\mathfrak{m})$. 
Now, $\mathfrak{m}$ is a left ideal of the super algebra $A$ if and only if $\mathfrak{m} = \Sigma(\mathfrak{m})$. 
This last equality tells us that $M$ is in fact a module over the super algebra $A$, and \textit{a fortiori} $M \oplus M^{\Sigma}$ is also a module over $A$. 
If $\mathfrak{m} \neq \Sigma(\mathfrak{m})$, then $\mathfrak{m} + \Sigma(\mathfrak{m}) = A$ and the Chinese Remained Theorem 
(\textit{cf.}~\cite{AF92}, Exercise 6.18) gives us an isomorphism 
\[     M \oplus M^{\Sigma} \simeq A/\mathfrak{m} \oplus A/\Sigma(\mathfrak{m}) \simeq A/(\mathfrak{m} \cap \Sigma(\mathfrak{m}))     \]
of modules over the algebra $\O(A)$. 
Since $\mathfrak{m} \cap \Sigma(\mathfrak{m})$ is a left ideal of the super algebra $A$, the claim holds. 

We now prove that $M_{+}$ and $M_{-}$ are simple over $\O(\U(\g))$. 
It is clear that $M_{+}^{\Sigma} \simeq M_{-}$, so the direct sum $M_{+} \oplus M_{-}$ has the structure of a module over the super algebra $\U(\g)$. 
Moreover, by construction, we have the isomorphisms $M_{+} \simeq \U(\g)/(\U(\g) J^{+})$ and $M_{-} \simeq \U(\g)/(\U(\g) J^{-})$ of modules over $\O(\U(\g))$. 
On the other hand, it can be easily verified that there are isomorphisms $\U(\g) J^{-}/(\U(\g) J) \simeq \U(\g)/(\U(\g) J^{+})$ and 
$\U(\g) J^{+}/(\U(\g) J) \simeq \U(\g)/(\U(\g) J^{-})$ of modules over the algebra $O(\U(\g))$, which in turn implies the isomorphism 
$M_{+} \oplus M_{-} \simeq \U(\g)/(\U(\g) J)$ of modules over $\O(\U(\g))$. 
Hence, we shall identify $M_{+}$ and $M_{-}$ with their images inside of $\U(\g)/(\U(\g) J)$ given by the previous isomorphisms. 
Furthermore, the isomorphism $\Sigma$ on $\U(\g)$ induces an isomorphism on $\U(\g)/(\U(\g) J)$ which sends $M_{+}$ to $M_{-}$. 
Since $J$ is a maximal left ideal of the super algebra $\U(\g)$, $\U(\g)/(\U(\g) J)$ is a simple module over $\U(\g)$. 
Let us choose $N$ a simple submodule of $M_{+}$. 
Then $\Sigma(N) \simeq N^{\Sigma}$ is a simple submodule of $M_{-} = \Sigma(M_{+})$ and by the previous paragraph $N \oplus \Sigma(N)$ has the structure 
of a module over the super algebra $\U(\g)$. 
Since $\U(\g)/(\U(\g) J)$ is a simple module over $\U(\g)$, then it coincides with $N \oplus \Sigma(N)$, which in turn implies that 
$M_{+} = N$ and $M_{-} = \Sigma(N)$. 
Hence, they are simple modules over $\O(\U(\g))$ and the proposition is thus proved. 
\qed

\subsubsection{Some results on stabilizers}

We want to describe the stabilizers of the primitive ideals of $\U(\g)$. 
In order to do that, we first consider the following simple result. 
\begin{lema}
\label{lema:6.2.6}
Let $\g$ be a nilpotent super Lie algebra, $\k$ an ideal of $\g$ such that there exists a homogeneous element $x \in \g$ satisfying that 
$\g$ is generated by $\k$ and $x$, and $\lambda \in \L_{\g}$ a functional satisfying that $\lambda(x) = \lambda([x,\k])=0$. 
Denote $\lambda' \in \L_{\k}$ the restriction of $\lambda$ to $\k$. 
Then $\st(I(\lambda'),\g)=\g$. 
\end{lema}
\noindent\textbf{Proof.}
The proof given in \cite{Dix96}, Lemma 6.2.6, (iii), also holds in this case, taking into account that there exists a polarization invariant under the action of $\g$ 
by Remark \ref{rem:polinv}, (i), and that $I(\lambda')$ is invariant under $\g$ if and only if it is invariant under $\ad(x)$. 
\qed

As a direct consequence from the previous lemma we obtain that (\textit{cf.}~\cite{Dix96}, Prop. 6.2.8):
\begin{prop}
\label{prop:6.2.8}
Let $\g$ be a nilpotent super Lie algebra, $\k$ an ideal of $\g$ and $\lambda \in \L_{\g}$. 
Denote by $\g'$ the super vector space formed by the elements $x \in \g$ satisfying that $\lambda([x,\h]) = 0$. 
Then, $\st(I(\lambda),\g) \supseteq \g' + \k$. 
\end{prop}
\noindent\textbf{Proof.}
It is direct that $\k \subset \st(I(\lambda),\g)$, so let us prove that $\g' \subset \st(I(\lambda),\g)$. 
Take $x \in \g'$ and consider the subalgebra $\h$ of $\g$ generated by $\k$ and $x$, which is nilpotent since $\g$ is. 
Lemma \ref{lema:6.2.6} tells us that $x \in \st(I(\lambda),\g)$ and the proposition is proved.
\qed



\noindent Estanislao Herscovich,
\\
Fakult\"at f\"ur Mathematik,
\\
Universit\"at Bielefeld,
\\
D-33615 Bielefeld,
\\
Germany,
\\
\href{mailto:eherscov@math.uni-bielefeld.de}{eherscov@math.uni-bielefeld.de}

\end{document}